\documentclass[11pt,reqno]{amsart}
\usepackage{microtype}
\usepackage{xcolor}
\usepackage[shortlabels]{enumitem}
\setlist{nolistsep}

\newtheorem{thm}{Theorem}[section]

\newtheorem{cor}[thm]{Corollary}
\newtheorem{prop}[thm]{Proposition}
\theoremstyle{definition}
\newtheorem{defn}[thm]{Definition}
\newtheorem{case}{Case}[thm]

\newtheorem{note}{Claim}[thm]
\newtheorem{rem}[thm]{Remark}
\newtheorem{exmp}[thm]{Example}
\numberwithin{equation}{section}

\title[Strong Numerical radius orthogonality]{Strong Numerical radius orthogonality  on bounded linear operators}
\author{Shaoxing Sun, Xiaofei Qi, Jinchuan Hou}
\address
[Sun]{
  Fundamentals Department, Air Force Engineering University, Xi'an, 710051, P. R. China}

\address[Qi]{
  School of Mathematics and Statistics, Shanxi University, Taiyuan 030006,
  P. R. China; Key Laboratory of Complex Systems and Data Science of Ministry of Education,
	Shanxi University, Taiyuan  030006,  Shanxi, China}

\address[Hou]{
  College of Mathematics, Taiyuan University of Technology, Taiyuan 030024,
P.R. China}

\email{ssx15603406512@163.com;
  xiaofeiqisxu@aliyun.com;
  jinchuanhou@aliyun.com}

\thanks{This work is partially supported by National Natural Science Foundation of
  China (12571138, 12071336), Shanxi Scholarship Council of China (2025-001)  and the Natural Science Basic Research Program of Shaanxi Province  (S2025-JC-QN-1854).}

\begin{document}

\begin{abstract}
 Orthogonality is an important tool for studying  bounded linear  operators. In this paper, we  introduce the concept of strong numerical radius orthogonality of operators and then, discuss the properties about numerical radius orthogonality, strong numerical radius orthogonality of operators as well as relationship between them.
\end{abstract}

\maketitle

\section{Introduction}

Orthogonality is a fundamental concept in Hilbert spaces and even operator algebras.
There is a natural orthogonality relation in an inner space $H$ with inner product $\langle \cdot,\cdot\rangle$. For any $x,y\in H$,  $x$ and $y$ are called orthogonal, denoted by $x\bot y$, if $\langle x,y\rangle =0$. Obviously, the orthogonality is symmetric, that is, $x\bot y$ if and only if $y\bot x$.
For any normed space   $X$ over $\mathbb{F} $ with norm $\|\cdot\|$, there exsits some kind of orthogonality, too.  For any $x,y\in X$,  $x$ is said to be Birkhoff-James orthogonal (B-J orthogonal) to $y$, denoted by $x\bot _B y$, if $\|x+\lambda y \| \geq \| x\| $ holds for all $\lambda \in\mathbb{F} $ \cite{B,J}. In general, B-J orthogonality is asymmetric, which is important in determining the geometry of the space. Also note that the above two kinds of orthogonality are equivalent in the case that $X$ is an inner space. The B-J orthogonality enables the concept of orthogonality to be applied to broader spaces,  and is of great significance for the study of the geometric properties and topological structures of the normed spaces \cite{GKSZ,MMQR,T}. It is worth noting that the series of studies on B-J orthogonality of operators  initiated  from the work of R. Bhatia and P. $\breve{S}$emrl \cite{BS}. After then, this topic were  widely discussed 
\cite{BG,PSG,PMW,T2}.

In \cite{AR}, a stronger version of B-J orthogonality was introduced in Hilbert C$^*$-modules.  Assume that $\mathcal M$ is a Hilbert
$\mathcal{A} $-module with $\mathcal A$ a C$^*$-algebra. Recall that, for  $x,y\in\mathcal M$, $x$ is strong B-J orthogonal to $y$, denote by $x\perp ^s_By$,  if $\|x+ya \|\geq \| x\| $ holds for all $a\in\mathcal{A} $. Note that
$\mathcal{B} (H)$  (the algebra of all bounded linear operators on a complex Hilbert space $H$) is  a Hilbert $\mathcal{B} (H)$-module. Thus, for  $A,B\in\mathcal{B} (H)$, we say that $A\perp ^s_BB$ if $\| A+BC\| \geq  \| A\| $ holds for all $C\in\mathcal{B} (H)$. The addition of multiplication can also facilitate a better study of the algebraic structure of the related operator algebras. For more informations about strong B-J orthogonality, the readers can see \cite{AR2,KST} and the references therein.

On the other hand, the numerical range and numerical radius are useful tools in studying matrices and operators, and have important  applications in many research areas \cite{W2}.
Recall that the numerical range and the numerical radius of $A\in \mathcal{B} (H)$ are respectively
$W(A)=\{\langle Ax,x\rangle :x\in H, \ \| x\| =1 \}$ and  $w(A)=\sup\{|\lambda |:\lambda \in W(A)\}$. Note that $w(\cdot)$ defines a norm on $\mathcal{B} (H)$, which is equivalent to the operator norm $\|\cdot\|$. Recently, the authors in \cite{MPS} introduced the definition of numerical radius orthogonality on $\mathcal B(H)$, that is,  $A\in\mathcal{B} (H)$ is said to be numerical radius orthogonal to $B\in\mathcal{B} (H)$, denoted by $A\bot _wB$, if $w(A+\lambda B)\geq w(A)$ holds for all $\lambda \in\mathbb{C} $.  It was  proved that $\bot _B$ and $\bot _w$ are not comparable;  in addition, some useful properties about numerical radius orthogonality were obtained \cite{MPS}. 

\begin{prop}\label{1.1} {\rm (\cite{MPS})}
	Assume that $H$ is a complex Hilbert space and
 $T,A\in\mathcal{B} (H)$. The following statements are equivalent.

 {\rm (1)} $T\perp _wA$.

{\rm (2)} $T^*\perp _wA^*$.

 {\rm (3)} $\alpha T\perp _w\beta A$ for all $\alpha ,\beta \in\mathbb{C} \setminus \{0\}$.

{\rm (4)} For each $\theta \in[0,2\pi)$, there exists a sequence $\{x_n^\theta \}\subseteq {\rm ball}H=\{x\in H:\| x\|=1 \}$ such that $\lim\limits _{n\rightarrow \infty }|\langle Tx_n^\theta ,x_n^\theta \rangle |=w(T)$ and $\lim\limits _{n\rightarrow \infty }{\rm Re}\{e^{-i\theta }\langle Tx_n^\theta ,x_n^\theta \rangle\overline{\langle Ax_n^\theta ,x_n^\theta \rangle}  \}\geq 0$.
\end{prop}

After then, a large number of articles began to study related topics of numerical radius orthogonality,  
for exmaple, see \cite{AK,BFP,GMPS,M,RS,TAK} and the references therein. These studies not only deepen the understanding of the geometric properties of the operators, but also refine some inequalities of the numerical radius.

Inspired by the above works, in this paper,  we introduce the concept of strong numerical radius orthogonality.

\begin{defn}\label{de1}

Let $H$ be a complex Hilbert space with $\dim H\geq 2$ and $A,B\in\mathcal{B} (H)$. We say that $A$ is  strong numerical radius orthogonal to $B$, denoted by $A\bot^s_w B$, if $w(A+BC)\geq w(A)$ holds for all $C\in\mathcal{B} (H)$.

\end{defn}

\begin{rem}\label{rem1}
     Note that, if we define $A\bot^{s}_wB$ by $w(A+CB)\geq w(A)$ for any $C\in\mathcal{B} (H)$, then it is easily seen that  $w(A^*+B^*C^*)\geq w(A^*)$ holds. So, without loss of generality, we adopt the definition of strong numerical radius orthogonality in Definition \ref{de1}.
    \end{rem}

The purpose of this  paper is to discuss the properties of numerical
radius orthogonality and strong numerical
radius orthogonality on bounded linear operators. The paper is organized as follows.
In Sections 2-3, we give useful properties of numerical radius orthogonality and strong numerical radius orthogonality of operators, respectively. In Section 4,
we give some characterizations of special bounded linear operators by the strong numerical radius orthogonality. 

We close this section by fixing some notions and notations.  Assume that  $H$ is a complex Hilbert space with $\dim H\geq 2$. Denote by $\mathcal{B} (H)$ and $\mathcal{B} _s(H)$ the algebra of all bounded linear operators  on  $H$ and the real  Jordan algebra of all self-adjoint operators in $\mathcal{B} (H)$, respectively. Let $\mathcal{M} _n(\mathbb{C} )$ be the algebra of all $n\times n$ complex matrices   and $\mathsf{H} _n$  the set of all Hermitian matrices in $\mathcal{M} _n(\mathbb{C} )$. For $A\in\mathcal{B} (H)$, $\sigma (A)$,   $\sigma_l(A)$, $\sigma_r(A)$,  $\sigma _p(A)$, $\sigma _{ap}(A)$,  $\ker A$, ${\rm ran}\,A$,  ${\rm rank}A$, $A^t$ and $\bar A$ respectively stand for the spectrum, the left spectrum, the right spectrum, the point spectrum, the approximate point spectrum, the kernel, the range, the rank,  the transpose and the conjugation  of $A$ with respect to an arbitrarily fixed orthonormal basis of $H$.   Note that $A^*={\bar{A}}^t$. Let ${\rm ball}\,H:=\{x\in H:\| x\|=1 \}$. For any subset $M\subset H$,  we will use $\overline{M} $, ${\rm co}\,M$,  ${\rm span}M$ and $M^\perp$ to denote the closure of $M$, the convex hull of  $M$, the linear span of the vectors in $M$ and the orthogonal complement subpace of $M$, respectively. Particularly, if $M=\{x\}$, write span$\{x\}=[x]$.

\section{Numerical radius orthogonality of bounded linear operators}

In this section, we discuss further the properties on numerical radius orthogonality of operators.

It was shown in  \cite[Theorem 2.2]{T2} that, for any $A\in{\mathcal B}(H)$,    $I\perp_BA$ if and only if $0\in\overline{W(A)} $. The next proposition indicates that this conclusion also holds for numerical radius orthogonality.

\begin{prop}\label{prop1}
Assume that   $A\in\mathcal{B} (H)$. Then  $I\perp _wA$ if and only if $0\in \overline{W(A)} $.
\end{prop}

\begin{proof}
Assume that $A\in{\mathcal B}(H)$. If $I\perp _wA$, by  Proposition \ref{1.1},
for any $\theta \in[0,2\pi)$, there exists a sequence $\{x_n^\theta \}\subset {\rm ball}\,H$ such that
    \begin{equation}\label{2.1}
       \lim_{n \to \infty}  |\langle x_n^\theta ,x_n^\theta \rangle |=1 \mbox{ and } \lim_{n \to \infty}{\rm Re}\,\{e^{-i\theta }\overline{\langle Ax_n^\theta,x_n^\theta\rangle } \}\geq 0.
      \end{equation}
Suppose on the contrary that $0\notin \overline{W(A)} $, and then there exist some $\theta_0 \in(-\frac{\pi}{2},\frac{\pi}{2})$ and $\varepsilon \in\mathbb{R}\backslash \{0\}$ such that either
    \begin{equation}\label{2.2}
        W(A)\subseteq \{x+iy:y>(\tan \theta _0)x+\varepsilon \},\quad  0\in\{x+iy:y<(\tan \theta _0)x+\varepsilon \};
    \end{equation}
    or
  \begin{equation}\label{2.2-2}
  	W(A)\subseteq \{x+iy:y<(\tan \theta _0)x+\varepsilon \},\quad  0\in\{x+iy:y>(\tan \theta _0)x+\varepsilon \}. \end{equation}
If Eq.(\ref{2.2}) holds, by taking
    $ \theta _1=
        \frac{\pi}{2}-\theta _0$ if $\theta _0\in[0,\frac{\pi}{2})$ and
       $ \theta _1=\pi+\theta _0$ if $ \theta _0\in(-\frac{\pi}{2},0)$, one achieves
    $${\rm Re}(e^{-i\theta_1} \overline{\langle Ax,x\rangle})={\rm Re }(e^{i\theta_1 }\langle Ax,x\rangle) <-\frac{|\varepsilon |}{\sqrt{1+\tan^2 \theta _0}}<0 {\rm \ for \ all \ } x\in {\rm ball}\,H.$$
This  contradicts Eq.(\ref{2.1}).  If Eq.\eqref{2.2-2} holds, a simiar argument can also get a contradiction. So $0\in \overline{W(A)} $.

Conversely, assume $0\in\overline{W(A)} $; then there exists a sequence $\{x_n\}\subset {\rm ball}\,H$ with $\|x_n\|=1$ such that $\lim\limits _{n\to \infty }\langle Ax_n,x_n\rangle =0$. It follows that
    $$\lim _{n\rightarrow \infty }|\langle x_n,x_n\rangle |=1=w(I) \ {\rm and } \ \lim _{n\rightarrow \infty }{\rm Re}(e^{i\theta }\langle Ax_n,x_n\rangle )=0.$$
By Proposition \ref{1.1},  we have $I\perp _wA$. The proof is completed.
\end{proof}

Note that the numerical radius orthogonality is asymmetric, that is, for  $A,B\in\mathcal{B} (H)$, $A\bot _wB$ cannot imply  $B\bot _wA$ in general.  Based on this point, the authors in  \cite{GMPS} introduced the concept of ``nr-symmetric",  where an operator $T\in\mathcal{B} (H)$ is called nr-right symmetric if,  for any $A\in\mathcal{B} (H)$, $A\bot _wT$ implies $ T\bot _wA$. It was proved  in \cite{GMPS} that, if $2\leq\dim H<\infty $ and $B\in\mathcal{B} _s(H)$ is  nr-right symmetric, then $B\in\mathbb{R} I$.  Here, we show that this result also holds for the infinite-dimensional case. In fact, we slightly generalize the result in \cite{GMPS} by relaxing any $A\in{\mathcal B}(H)$ to any $A\in{\mathcal B}_s(H)$  for the case $\dim H\geq 3$.

\begin{prop}\label{prop2}
Assume that $\dim H\geq 3$ and $B\in\mathcal{B} _s(H)\backslash \{0\}$.  The following statements are equivalent.

{\rm (i)} $B$ is nr-right symmetric, that is, for any  $A\in\mathcal{B}(H)$, $A\perp _wB$ implies $B\perp _wA$.

{\rm (ii)} For any $A\in\mathcal{B} _s(H)\backslash \{0\}$, $A\perp _wB$ implies $B\perp _wA$.

{\rm (iii)}  $B=\lambda I$ with $\lambda \in\mathbb{R} \backslash \{0\}$, where $\mathbb R$ is the real number field.

\end{prop}

\begin{proof}
   (i)$\Rightarrow $(ii): Obvious.

   (iii)$\Rightarrow $(i): Note that, by \cite[Proposition 2.4]{TAK},  for any $A\in\mathcal{B} (H)$,  $A\perp _wI$ implies $ I\perp _wA$. Combining this fact and Proposition \ref{1.1},  (iii)$\Rightarrow $(i) is clearly true.

(ii)$\Rightarrow $(iii): Since $B$ is self-adjoint, we can assume $\overline{W(B)} =[\lambda _1,\lambda _2]$ with $\lambda _1,\lambda _2\in\mathbb{R} $ and $\lambda _1\leq \lambda _2$.  Write  $B=\int z{\rm d}E(z)$ the spectral decomposition.

We first prove the following claim.

{\bf Claim.}
     $0\notin\overline{W(B)} $.

On the contrary, assume  $0\in\overline{W(B)} $.
By Proposition \ref{prop1},  we have $I\perp _wB$, which and (ii) give $B\perp _wI$. Then, by Proposition \ref{1.1}, for any $\theta \in[0,2\pi)$, there exists a sequence $\{y_n^\theta \}\subset {\rm ball}\,H$ such that
$$\lim _{n\rightarrow \infty }|\langle By_n^\theta,y_n^\theta\rangle |=w(B) {\rm \ and \ }\lim _{n\rightarrow \infty }{\rm Re}(e^{-i\theta }\langle By_n^\theta,y_n^\theta\rangle )\geq 0.$$
It follows from the arbitrariness of $\theta$ that $\lambda _1=-\lambda _2$. By the equivalence of (1) and (3) in Proposition \ref{1.1}, we may assume that $\lambda _2=1$, and thus $\overline{W(B)} =[-1,1]$.

{\bf Case 1.}
   $\sigma (B)\backslash \{-1,1\}\neq\emptyset $.

In this case,
take $\varepsilon _0> 0$ such that $-1<-1+\varepsilon _0\leq 0\leq 1-\varepsilon _0<1$.  Decompose $H$ into four subspaces:
$H=H_1\oplus H_2\oplus H_3\oplus H_4$, where
   $$H_1=E([-1,-1+\varepsilon _0))H, H_2=E([-1+\varepsilon _0,0))H, H_3=E([0,1-\varepsilon _0))H, H_4=E([1-\varepsilon _0,1))H.$$
We may assume $\dim H_3,\dim H_4>0$; otherwise, replace $B$ with $-B$ without loss of generality. Thus, under the above space decomposition, $B$ has the form
   $$B=B_1\oplus B_2\oplus B_3\oplus B_4.$$

Take
   $$A=I_1\oplus I_2\oplus I_3\oplus (-I_4).$$
Considering  $\{x_n=0\oplus 0\oplus x_n^3\oplus 0 \}$ with $\{x_n^3\}\subset{\rm ball}\,H_3$ for $\theta \in[0,\frac{\pi}{2})\cup (\frac{3\pi}{2},2\pi)$ and  $\{x_n=0\oplus 0\oplus 0\oplus x_n^4 \}$ with $\{x_n^4\}\subset {\rm ball}\,H_4$ for $\theta \in[\frac{\pi}{2},\frac{3\pi}{2}]$, we have
   $$\lim_{n \to \infty} |\langle Ax_n,x_n\rangle | =1=w(A){\rm \ and \ }\lim_{n \to \infty}{\rm Re}\{e^{-i\theta }\langle Ax_n,x_n\rangle \langle Bx_n,x_n\rangle \}\geq 0.  $$
   It follows that $A\bot _wB$.

However,  $B\not\perp_wA$.  In fact, for any $\{y_n\}\subset  {\rm ball}\,H$ with $\lim\limits _{n\rightarrow \infty }|\langle By_n,y_n\rangle |=w(B)=1$, write
$$y_n=\alpha _nu_n\oplus \beta _nv_n\oplus \gamma _nz_n,$$ where $u_n\in {\rm ball}\,H_1$, $v_n\in {\rm ball}\,(H_2\oplus H_3)$, $z_n\in {\rm ball}\,H_4$, $\alpha _n,\beta _n,\gamma _n\in\mathbb{C} $ and $|\alpha _n|^2+|\beta _n|^2+|\gamma _n|^2=1$. Then
   \begin{equation}\label{2.3}
    \langle By_n,y_n\rangle =|\alpha _n|^2\langle B_1u_n,u_n\rangle +|\beta _n|^2\langle (B_2\oplus B_3)v_n,v_n\rangle+|\gamma _n|^2\langle B_4z_n,z_n\rangle ,
   \end{equation}
   which implies
   $$|\langle By_n,y_n\rangle |\leq |\alpha _n|^2+|\beta _n|^2(1-\varepsilon _0)+|\gamma _n|^2=1-\varepsilon _0|\beta _n|^2.$$
 So $1\leq 1-\varepsilon _0\lim\limits_{n \to \infty}  |\beta _n|^2$, and hence $\lim\limits_{n \to \infty}  |\beta _n|=0$. As $\{\langle By_n,y_n\rangle \}$ is a sequentially compact set, there exists a convergent subsequence.  Without loss of generality, still denote by $\{\langle By_n,y_n\rangle \}$. Thus, we conclude that either $\lim\limits _{n\rightarrow \infty }\langle By_n,y_n\rangle=1$ or $\lim\limits _{n\rightarrow \infty }\langle By_n,y_n\rangle=-1$.
 If $\lim\limits _{n\rightarrow \infty }\langle By_n,y_n\rangle=1$, by Eq.(\ref{2.3}), one has
$$\begin{array}{rl}1=& \lim\limits_{n \to \infty}  \langle By_n,y_n\rangle= \lim\limits_{n \to \infty} (|\alpha _n|^2\langle B_1u_n,u_n\rangle +|\gamma _n|^2\langle B_4z_n,z_n\rangle)
  \\
=&\lim\limits_{n \to \infty}((|\gamma _n|^2-1)\langle -B_1u_n,u_n\rangle+|\gamma _n|^2\langle B_4z_n,z_n\rangle )\\
\leq &\lim\limits_{n \to \infty}2|\gamma _n|^2-1\leq 1,
      \end{array}$$
which entails $\lim\limits_{n \to \infty}  |\gamma  _n|=1$. If $\lim\limits _{n\rightarrow \infty }\langle By_n,y_n\rangle=-1$, an analogous discussion leads to $\lim\limits_{n \to \infty}  |\alpha  _n|=1$. Now, by taking $\theta \in(0,\frac{\pi}{2})$, one always  has
   $$\lim _{n\rightarrow \infty }{\rm Re}\{e^{-i\theta }\langle By_n,y_n\rangle \langle Ay_n,y_n\rangle \}=-\cos \theta <0.$$
It follows from Proposition \ref{1.1}  that $B\not\perp_wA$, contradicting to (ii).

{\bf Case 2.}
  $\sigma (B)\backslash \{-1,1\}=\emptyset  $.

   In this case, $-1,1$ are isolated points of $\sigma (B)$,  and so  $-1,1\in\sigma _p(B)$. Let $$H_1=\ker (I-B)\ \ {\rm and}\ \ H_2=\ker (-I-B).$$ Still, as   $\dim H\geq 3$, without loss of generality, we may assume $\dim H_1\geq 2$. Take vectors $x,y\in {\rm ball} \ H_1$ and $z\in {\rm ball} \ H_2$ with $\langle x,y\rangle =0$. With respect to the space decomposition $H={\rm span}\{x,y,z\}\oplus (H_1\ominus {\rm span}\{x,y\})\oplus  (H_2\ominus [z])$, $B$ can be written as
$$B=B_1\oplus I_1\oplus (-I_2) {\rm \ with \ }B_1=\left(
  \begin{array}{ccc}
      1 & 0 & 0  \\
      0 & 1 & 0  \\
      0 & 0 & -1
    \end{array}
  \right).$$

Take
$$A=A_1\oplus (-\frac{1}{2}I_1)\oplus \frac{1}{2}I_2,$$
where $A_1=u\otimes u-\frac{2}{3}v\otimes v$ with $u=\frac{1}{\sqrt{6}}x+\frac{1}{\sqrt{3}}y+\frac{1}{\sqrt{2}}z$ and $v=\frac{1-\sqrt{14}}{3\sqrt{3}}x+\frac{2+\sqrt{14}}{3\sqrt{6}}y-\frac{1}{3}z$. Note that $\langle u,v\rangle =0$ and $\| u\|= \|v \|=1 $. A simple calculation yields $w(A)=1=\langle Au,u\rangle$, and then
$${\rm Re}\{e^{-i\theta }\langle Au,u\rangle \langle Bu,u\rangle \}=0{\rm \ holds \ for \ all  \ } \theta \in[0,2\pi).$$
It follows from Proposition \ref{1.1} that $A\perp _wB$.

On the other hand, note that
$$A_1=\left(
  \begin{array}{ccc}
      \frac{1}{6}-\frac{2}{3}\cdot \frac{15-2\sqrt{14}}{27}                & \frac{1}{3\sqrt{2}}+\frac{2}{3}\cdot \frac{12+\sqrt{14}}{27\sqrt{2}} & *                        \\
      \frac{1}{3\sqrt{2}}+\frac{2}{3}\cdot \frac{12+\sqrt{14}}{27\sqrt{2}} & \frac{1}{3}-\frac{2}{3}\cdot \frac{9+2\sqrt{14}}{27}                 & *                        \\
      *                                                                    & *                                                                    & \frac{1}{2}-\frac{2}{27}
    \end{array}
  \right).$$
This means $\langle Ax',x'\rangle <0 $ for all unit vectors $x'\in{\rm span}\{x,y\}$. For any $q\in H_1$, wirte $q=q_1\oplus q_2$ with $q_1\in{\rm span}\{x,y\}$ and $q_2\in H_1\ominus {\rm span}\{x,y\}$. As ${\rm span}\{x,y,z\}$ reduces $A$, we have
\begin{equation}\label{2.4.1}
  \langle Aq,q\rangle =\langle A_1q_1+(-\frac{1}{2}I_1)q_2,q_1+q_2\rangle=\langle A_1q_1,q_1\rangle-\frac{1}{2}<0{\rm \ for \ any \ }q\in H_1.
\end{equation}
Similarly, one can get
{\begin{equation}\label{2.4}
\langle Aq,q\rangle >0 {\rm \ for \ any \ } q\in H_2.
\end{equation}}
Now, take a sequence $\{y_n=\alpha _nu_n\oplus \beta _nv_n\}\subset  {\rm ball}\,H$ with $u_n\in {\rm ball}\,H_1$, $v_n\in {\rm ball}\,H_2$, $\alpha _n,\beta _n\in\mathbb{C} $ and $|\alpha _n|^2+|\beta _n|^2=1$, such that $\lim\limits _{n\rightarrow \infty }|\langle By_n,y_n\rangle |=w(B)=1$. Then
$$\lim\limits _{n\rightarrow \infty }|\langle By_n,y_n\rangle| =\lim\limits _{n\rightarrow \infty }(||\beta _n|^2-|\alpha _n|^2|)=1.$$
So either $\lim\limits_{n \to \infty}  |\alpha _n|=0$ or $\lim\limits_{n \to \infty}  |\beta _n|=0$. By taking $\theta \in(0,\frac{\pi}{2})$ and using Eqs.(\ref{2.4.1})-(\ref{2.4}), we can infer that
$$\lim _{n\rightarrow \infty }{\rm Re}\{e^{-i\theta }\langle By_n,y_n\rangle \langle Ay_n,y_n\rangle \}<0.$$
This implies  $B\not\perp _wA$,  a  contradiction again.

Hence the claim holds.

Now, if $B\notin \mathbb{R} I$, then $\lambda _1< \lambda _2$. Since $W(B)$ is a convex set, by Claim, we have either $[\lambda _1,\lambda _2]\subset (0,+\infty )$ or $[\lambda _1,\lambda _2]\subset (-\infty,0 )$. Without loss of generality,  assume $[\lambda _1,\lambda _2]\subset (0,+\infty )$.  Take $\varepsilon _0,\varepsilon _0'>0$ such that  $\lambda _1<\lambda _1+\varepsilon _0<\lambda _2-\varepsilon _0'<\lambda _2$, and let
$$H_1=E([\lambda _1,\lambda _1+\varepsilon _0))H,\ \ H_2=E([\lambda _1+\varepsilon _0,\lambda _2-\varepsilon _0'))H,\ \ H_3=E([\lambda _2-\varepsilon _0',\lambda _2])H.$$ Then, under the space decomposition $H=H_1\oplus H_2\oplus H_3$, we have
$B=B_1\oplus B_2\oplus B_3. $

Take $A=I_1\oplus 0\oplus (-I_3).$
It is clear that, for any $\theta \in[0,\frac{\pi}{2}]\cup [\frac{3\pi}{2},2\pi)$ and any $\{x_n=x_{n}^1\oplus 0\oplus 0\}$ with $\{x_{n}^1\}\subset {\rm ball}\, H_1$,  one has
$$\lim _{n\rightarrow \infty }|\langle Ax_n,x_n\rangle |=1=w(A) $$ and $$ \lim _{n\to \infty }{\rm Re }\{e^{-i\theta }\langle Ax_n,x_n\rangle \langle Bx_n,x_n\rangle \}=\lim _{n\rightarrow \infty }\cos \theta \langle Bx_n,x_n\rangle \geq 0;$$
and for any $\theta \in(\frac{\pi}{2},\frac{3\pi}{2})$ and   $\{x_n=0\oplus 0\oplus x_n^3\}$ with $\{x_n^3\}\subset{\rm ball}\,H_3$, one gets
$$\lim _{n\rightarrow \infty }|\langle Ax_n,x_n\rangle |=1=w(A)$$
and
$$ \lim _{n\to \infty }{\rm Re }\{e^{-i\theta }\langle Ax_n,x_n\rangle \langle Bx_n,x_n\rangle \}=\lim _{n\rightarrow \infty }-\cos \theta \langle Bx_n,x_n\rangle \geq 0.$$
Thus we have $A\perp _wB$.

On the other hand, for any $\{x_n\}\subset  {\rm ball}\,H$ with $\lim\limits _{n\rightarrow \infty }|\langle Bx_n,x_n\rangle |=w(B)=\lambda _2$, write
 $$x_n=\alpha _nu_n\oplus \beta _nv_n,$$where $u_n\in {\rm ball}\,(H_1\oplus H_2),\ v_n\in {\rm ball}\,H_3,\ \alpha _n,\beta _n\in\mathbb{C} \ {\rm with}\ |\alpha _n|^2+|\beta _n|^2=1.$ A direct calculation gives
$$\langle Bx_n,x_n\rangle =|\alpha _n|^2\langle (B_1\oplus B_2)u_n,u_n\rangle +|\beta _n|^2\langle B_3v_n,v_n\rangle \leq |\alpha _n|^2(\lambda _2-\varepsilon _0')+|\beta _n|^2\lambda _2=\lambda _2-|\alpha _n|^2\varepsilon _0',$$
which implies $\lambda _2\leq \lambda _2-\lim\limits _{n\rightarrow \infty }|\alpha _n|^2\varepsilon _0'$. So  $\lim \limits_{n\rightarrow \infty }|\alpha _n|^2=0$, and hence
$$\lim\limits _{n\rightarrow \infty }{\rm Re}\{e^{-i\theta }\langle Bx_n,x_n\rangle \langle Ax_n,x_n\rangle \}=-\lambda _2\cos \theta .$$
 It follows that,  for any $\theta \in(0,\frac{\pi}{2})$, and any  $\{x_n\}\subset  {\rm ball}\,H$ with $\lim\limits _{n\rightarrow \infty }|\langle Bx_n,x_n\rangle |=w(B)=\lambda _2$, we always have   $\lim\limits _{n\rightarrow \infty }{\rm Re}\{e^{-i\theta }\langle Bx_n,x_n\rangle \langle Ax_n,x_n\rangle \}<0$.    This means that $B\not\perp _wA$, again a contradiction to (ii).

Therefore, we must have $B\in\mathbb{R} I$, as desired.
\end{proof}

Note that, by   \cite[Theorem 3.5]{GMPS} and \cite[Proposition 2.4]{TAK},  the statements (i) and (iii) in Proposition \ref{prop2} is equivalent when $\dim H =2$. However, the following proposition demonstrates that the statement (ii) is not equivalent to (i) and (iii) in the case $\dim H =2$.  Thus Proposition \ref{prop2} is not true for the case $\dim H=2$.

\begin{prop}\label{prop2.1}
	 Let
	$B\in\mathsf{H} _2$. Then the following statements are equivalent.

{\rm (1)}  For any $A\in\mathsf{H} _2$, $A\perp _wB$ implies $B\perp _wA$.

 {\rm (2)}  Either $B\in\mathbb R I$ or $B=\lambda \left(
 \begin{array}{ccc}
 	1 & 0  \\
 	0 & -1
 \end{array}
 \right)$ with $\lambda\in\mathbb R\setminus\{0\}$.

	\end{prop}
	
\begin{proof}   (2)$\Rightarrow$(1): If $B\in\mathbb R I$,  by \cite[Proposition 2.4]{TAK}, the statement (1) holds.
	
If $B=\lambda \left(
	\begin{array}{ccc}
		1 & 0  \\
		0 & -1
	\end{array}
	\right)$ with $0\not=\lambda\in\mathbb R $, by Proposition \ref{1.1}, we may assume $\lambda=1$.  For any $A\in\mathsf{H} _2$,  write $A=\left(
	\begin{array}{ccc}
		a_1 & a_2  \\
		\bar a_2 & a_3
	\end{array}
	\right)$. It is straightforward to verify that $\sigma (A)=\{\lambda _1,\lambda _2\}$ with $\lambda _1=
	\frac{a_1+a_3+\sqrt{(a_1-a_3)^2+4|a_2|^2}}{2}$ and  $\lambda _2=\frac{a_1+a_3-\sqrt{(a_1-a_3)^2+4|a_2|^2}}{2}$. Without loss of generality,  assume $a_1+a_3\geq 0$; otherwise, replace $A$ with $-A$. Then we have $w(A)=\lambda _1$.
	
	In the sequel, assume that $A\perp _wB$.
	
	We claim that either
	(i) $a_2\neq 0$ and $a_1+a_3>0$, or
		(ii) $a_2= 0$ and $a_1=a_3$.
	Otherwise, if $a_2\neq 0$ and $a_1+a_3=0$, then
$A=\left(
\begin{array}{ccc}
	a_1 & a_2  \\
	\bar a_2 & -a_1
\end{array}
\right)$ and $\lambda_2=-\lambda _1=-\sqrt{a_1^2+|a_2|^2}$. Note that, for any $x\in{\rm ball}\,\mathbb{C} ^2$, if $|\langle Ax,x\rangle |=w(A)$ holds, then $x$ has the form
$$\frac{e^{i\delta }}{\sqrt{|a_2|^2+(\lambda -a_1)^2}}\begin{pmatrix}
	a_2\\
	\lambda-a_1
\end{pmatrix}\ \ {\rm with} \ \ \delta \in[0,2\pi) \ {\rm and}\ \lambda=\lambda_i, \ i=1,2.$$A direct calculation gives
$\langle Ax,x\rangle\langle Bx,x\rangle>0$. Thus, there exists some   $\theta _0\in(\frac{\pi}{2},\frac{3\pi}{2})$ such that for any  $x\in{\rm ball}\,\mathbb{C} ^2$ with $|\langle Ax,x\rangle |=w(A)$,
${\rm Re}\,\{e^{-i\theta _0}\langle Ax,x\rangle \langle Bx,x\rangle \}<0$
always holds. So $A\not\perp _wB$ by Proposition  \ref{1.1}.  If  $a_2=0$ and $a_1\not=a_3$, it is also easy to check  $A\not\perp _wB$.

If (i) holds,  $A$ has two distinct eigenvalues. For any $x\in{\rm ball}\,\mathbb{C} ^2$, if   $|\langle Ax,x\rangle |=w(A)$,  $x$ must be  the form  $\frac{e^{i\delta }}{\sqrt{|a_2|^2+(\lambda_1 -a_1)^2} }\begin{pmatrix}
		a_2\\
		\lambda _1-a_1
	\end{pmatrix}$ with $\delta \in[0,2\pi)$. By Proposition \ref{1.1} (1) and (4), we can infer  $\langle Bx,x\rangle =0$, and so $|a_2|^2=(\lambda _1-a_1)^2$. This implies
\begin{equation*}
	(a_1-a_3)^2=(a_1-a_3)\sqrt{(a_1-a_3)^2+4|a_2|^2}.
	\end{equation*}
	If $a_1\neq a_3$, then the above equation gives $a_2=0$, a contradiction. Hence $a_1=a_3$. Now, by taking  $x=\begin{pmatrix}
		1\\
		0
	\end{pmatrix}$ for any  $\theta \in[0,\frac{\pi}{2})\cup [\frac{3\pi}{2},2\pi)$ and $x=\begin{pmatrix}
	0\\
	1
	\end{pmatrix}$ for any $\theta \in[\frac{\pi}{2},\frac{3\pi}{2})$,
one can achieve
	${\rm Re}\,\{e^{-i\theta }\langle Bx,x\rangle \langle Ax,x\rangle \}\geq 0$. It follows from Proposition \ref{1.1} (1) and (4) that  $B\perp _wA$.
	
	If (ii) holds, then $A=a_1 I$. It is obvious that  $B\perp _w a_1I$.
	
	Hence, the statment (1) holds.

(1)$\Rightarrow$(2):  Assume that, for any $A\in\mathsf{H} _2$, $A\perp _wB$ implies $B\perp _wA$. By the equivalence of (1) and (3) in Proposition \ref{1.1}, we may assume that $B=\left(
\begin{array}{ccc}
	1 & 0  \\
	0 & \mu
\end{array}
\right)$ with $\mu \in\mathbb{R} $ and $|\mu |\leq 1$.

 If $-1<\mu\leq 0$,  by taking $A=I$, one has  $I\perp _wB$ and $B\not\perp _w I$, a contradiciton; if $0<\mu<1$,  by taking $A=\left(
\begin{array}{ccc}
	1 & 0  \\
	0 & -1
\end{array}
\right)$, one sees $A\perp _wB$ and $B\not\perp _wA$, a contradiciton again.

So $\mu =\pm1$, and hence, the statement (2) is true.

The proof is completed.\end{proof}

The next result reveals  that the B-J orthogonality and the numerical radius orthogonality are equivalent on $\mathcal{B} _s(H)$.

\begin{prop}\label{prop3}
For any $A,B\in\mathcal{B} _s(H)\backslash \{0\}$, we have
$$A\perp_w B\Leftrightarrow A\perp _BB.$$
\end{prop}

\begin{proof}  Note that $w(T)\leq \|T\|$ holds for all $T\in{\mathcal B}(H)$, and the equality holds when $T$ is self-adjoint.
Now, assume that $A,B\in\mathcal{B} _s(H)\backslash \{0\}$.
If $A\perp_w B$, then
    $$\| A+\lambda B\| \geq w(A+\lambda B)\geq w(A)=\| A\| {\rm \ holds \ for \ all \ }\lambda \in \mathbb{C} ,$$
    which implies $A\bot _BB$.

If $A\bot _BB$, then
$  \| A+rB\| \geq \| A\| $ holds for all $r\in\mathbb{R}.
$
For any $\lambda \in\mathbb{C} $, write $\lambda =t+is$ with $t,s\in\mathbb{R} $. Take  $\{x_n\}\subset {\rm ball}\,H$ with $\lim\limits_{n \to \infty} |\langle (A+tB)x_n,x_n\rangle|=w(A+tB)$. We have
    $$\begin{array}{rl}w^2(A+\lambda B)  \geq& \lim\limits_{n \to \infty}  |\langle (A+\lambda B)x_n,x_n\rangle |^2\\
    	=&\lim\limits_{n\to \infty}  |\langle (A+tB)x_n,x_n\rangle +i\langle sBx_n,x_n\rangle |^2 \\
        =&\lim\limits_{n\to \infty} (|\langle (A+tB)x_n,x_n\rangle |^2+|\langle sBx_n,x_n\rangle |^2)\\
        \geq& \lim\limits_{n \to \infty} |\langle (A+tB)x_n,x_n\rangle |^2   \\
          =&w^2(A+tB)=\| A+tB\| ^2\geq \| A\| ^2=w^2(A),
      \end{array}$$
and so  $A\bot _wB$.

The proof is finished.
\end{proof}

\section {Strong numerical radius orthogonality of bounded linear operators}

In this section, we give some properties about strong numerical radius orthogonality of bounded linear operators.

The first proposition gives the relationship between the  numerical radius orthogonality and strong numerical radius orthogonality.

\begin{prop}\label{prop3.1}
 For any  $A,B\in\mathcal{B} (H)$,
  $A\bot _w^s B$ if and only if $ A\bot _w BC$  for   all $C\in\mathcal{B} (H).$
\end{prop}

\begin{proof} Assume that $A,B\in\mathcal{B} (H)$. If $A\bot _w^s B$, we see that
	$$w(A)\leq w(A+B\cdot \lambda C)=w(A+\lambda BC)$$ holds for all $C\in\mathcal{B} (H)$ and $\lambda \in\mathbb{C} $, which means $A\bot _w BC$ for all $C\in\mathcal{B} (H)$.

If $A\bot_w BC$ holds for all $C\in\mathcal{B} (H)$, then $w(A)\leq  w(A+BC)$, which gives $A\bot _w^s B$.
\end{proof}

By Proposition \ref{3.1}, we see that, for any operators $A,B\in\mathcal B(H)$,
$A\bot _w^s B$ implies $ A\bot _w B$, but the inverse does not hold in general.

In addition, by using Proposition \ref{3.1},  we can get a sufficient condition of $A\bot_w^s B$.

\begin{prop}\label{prop3.2}
 Assume that $A,B\in\mathcal{B} (H)$. If there exists a sequence $\{x_n\}\subset  {\rm ball}\,H$ such that $\lim\limits _{n\rightarrow \infty }|\langle Ax_n,x_n\rangle |=w(A)$ and $\lim\limits _{n\rightarrow \infty }\| B^*x_n\|=0$, then $A\bot_w^s B$.
    \end{prop}
    \begin{proof}
        Assume that $A,B\in\mathcal{B} (H)$ and there exists a sequence $\{x_n\}\subset  {\rm ball}\,H$ such that $\lim\limits _{n\rightarrow \infty }|\langle Ax_n,x_n\rangle |=w(A)$ and $\lim \limits_{n\rightarrow \infty }\| B^*x_n\| =0$. Then, for any $\theta \in[0,2\pi)$, one has
        $$\lim _{n\rightarrow \infty }{\rm Re}\{e^{-i\theta }\langle Ax_n,x_n\rangle \overline{\langle BCx_n,x_n\rangle } \}=0 {\rm  \ for \ all \ }C\in\mathcal{B} (H).$$
  This, together with Proposition \ref{1.1} and  Proposition \ref{3.1} gives  $A\perp _w^sB$.
    \end{proof}

    Particularly, in  Proposition \ref{prop3.2},  if  $A$ is a normaloid operator, then the sufficient condition  is also  necessary. Recall that an operator $A\in\mathcal{B} (H)$ is called normaloid if $\| A\| =w(A)$ \cite{W1}. Obviously,  every normal  operator must be normaloid.

  \begin{prop}\label{prop3.8}
  	 Assume that $A,B\in\mathcal{B}(H)$ with $0\not=A$   normaloid. Then the following statements are equivalent.
  	
  {\rm 	(i)}  $A\bot ^s_w B$.
  	
   {\rm 	(ii)}  There exist a sequence $\{x_n\}\subset {\rm ball}\,H$ such that $\lim\limits_{n \to \infty}|\langle Ax_n,x_n\rangle |=w(A) $ and $\lim  	\limits_{n \to \infty}\| B^*x_n\| =0$.

  \end{prop}

  \begin{proof}
  	(ii)$\Rightarrow $(i): It is obvious by Proposition \ref{prop3.2}.
  	
  	(i)$\Rightarrow $(ii): On the contrary, assume that for any sequence $\{x_n\}\subset {\rm ball}\,H$ with $\lim\limits_{n \to \infty}|\langle Ax_n,x_n\rangle |=w(A)$,   $\lim\limits_{n \to \infty}\| B^*x_n\| \neq 0$.
  Since $\{\langle Ax_n,x_n\rangle\}$ is a sequentially compact set, there exist a subsequence $\{\langle Ax_{n_i},x_{n_i}\rangle\}$ and $\lambda \in\mathbb{C} $ such that $\lim\limits_{n \to \infty}\langle Ax_{n_i},x_{n_i}\rangle=\lambda $.  So $|\lambda |=w(A)$.
    Let $Ax_{n_i}=\lambda _{n_i}x_{n_i}+\mu _{n_i}y_{n_i}$ with $\{y_{n_i}\}\subset{\rm ball}\,H$ and $\langle x_{n_i},y_{n_i}\rangle =0$, $\forall n$. Then $\lim\limits_{n \to \infty} \lambda _{n_i}=\lambda $. Note that $\| A\| =w(A)=|\lambda |$. So  $\lim\limits_{n \to \infty}|\mu _{n_i}| =0 $, and thus
  	$$\begin{array}{rl}
  		& \lim\limits_{n \to \infty} {\rm Re}\,\{e^{-i\pi}\langle Ax_{n_i},x_{n_i}\rangle\overline{\langle BB^*Ax_{n_i},x_{n_i}\rangle }  \} =\lim\limits_{n \to \infty} {\rm Re}\,\{e^{-i\pi}\lambda \bar\lambda _{n_i}\langle BB^*x_{n_i},x_{n_i}\rangle \}  \\
  		=& \lim\limits_{n \to \infty} {\rm Re}\,\{e^{-i\pi}|\lambda |^2\|B^*x_{n_i} \|^2 \}=-\lim\limits_{n \to \infty} |\lambda |^2\|B^*x_{n_i} \|^2<0.
  	\end{array}$$
It follows from  Proposition \ref{1.1} that $A\not\perp_w BB^*A$, and hence  $A\not\perp_w^sB$ by Proposition \ref{prop3.1}, a contradiction.

The proof is finished.
  \end{proof}

\if false
In fact, if both $A$ and $B$ can be choosen to be  self-adjoint operators, then Proposition \ref{prop3.8} still holds.\fi

\begin{cor}\label{prop3.6}
Assume that  $0\not=A\in\mathcal{B} (H)$ is either normal or  self-adjoint, and $B\in\mathcal B(H)$.  Then the following statements are equivalent.
	
 {\rm 	(i)}   $A\bot _w^sB$.
	
 {\rm 	(ii)}  There exists a sequence $\{x_n\}\subset  {\rm ball}\,H$ such that $\lim\limits _{n\rightarrow \infty }|\langle Ax_n,x_n\rangle |=w(A)$ and  $\lim\limits _{n\rightarrow \infty }\| B^*x_n\| =0$.
\end{cor}

Next, we consider  two kinds of special operators.

For $A\in\mathcal B(H)$, denote by $$M_A=\{x\in {\rm ball}\,H: |\langle Ax,x\rangle |=w(A)\},\ \ N_A=\{\lambda \in\mathbb C:\lambda \in W(A)\ {\rm and}\ |\lambda |=w(A)\},$$
$$N'_A=\{\lambda\in\mathbb C  :\lambda \in\overline{W(A)}\ {\rm and}\ |\lambda |=w(A) \}.$$
Let
$$\mathcal{E}=\{A\in\mathcal{B} (H):0\notin{\rm co}\,\{N'_A\}\}$$ and
$$\begin{array}{rl}\mathcal{F} =\{A\in\mathcal{M} _n(\mathbb{C} ): &{\rm either}\ 0\notin {\rm co}\,\{N_A\}; {\rm \ or \ for \ any \ subset \ }M_A'\subseteq M_A, \\
	&{\rm \ there \ exists \ a \ vector \ }x_0\in M_A' {\rm \ such \ that \ } x_0\notin {\rm span} \,\{M_A'\backslash \{x_0\} \}\},\end{array}$$
 where co$\{G\}$ stands for as usal the convex hull of the set $G$.


\begin{prop}\label{prop3.3}
	Assume that $A,B\in\mathcal{B} (H)$ and $A\perp _w^sB$. The following stamtents hold.
	
 {\rm 	(1)}   If $\dim H=\infty $ and $A\in\mathcal{E}$, then there exists a sequence $\{x_n\}\subset  {\rm ball}\,H$ such that $\lim\limits _{n\rightarrow \infty }|\langle Ax_n,x_n\rangle |=w(A)$ and $\lim\limits _{n\rightarrow \infty }\| B^*x_n\| =0$.
	
 {\rm 	(2)}  If $\dim H=n<\infty$ and  $A\in\mathcal{F} $, then there exists a vector $x\in {\rm ball}\,H$ such that $|\langle Ax,x\rangle |=w(A)$ and $B^*x=0$.
\end{prop}

\begin{proof}
	Assume that $A,B\in\mathcal{B} (H)$ and $A\perp _w^sB$.
	
	(1) 	If $\dim H=\infty $ and $A\in\mathcal{E} $, {  then there exist $\theta \in[0,2\pi)$ and $\varepsilon >0$ such that $$N'_A\subseteq \{w(A)e^{i\xi }: \xi \in[\theta +\varepsilon ,\theta +\pi-\varepsilon ]\}.$$} Suppose on the contrary that, for any $\{x_n\}\subset {\rm ball }\,H$ with $\lim\limits _{n\rightarrow \infty }|\langle Ax_n,x_n\rangle |=w(A)$, we have $\lim \limits_{n\rightarrow \infty }\|B^*x_n \|\neq 0 $.  By taking a subsequence, we may require that $\lim\limits _{n\rightarrow \infty }\langle Ax_n,x_n\rangle =e^{i\xi}w(A)$, which is an element in $N'_A$. Then, by taking   $\delta =\frac{\pi}{2}-\theta$, one achieves
	$$\lim _{n\rightarrow \infty }{\rm Re}\{e^{i\delta }\langle Ax_n,x_n\rangle \overline{\langle BB^*x_n,x_n\rangle}  \}=\lim _{n\rightarrow \infty }\|B^*x_n \|^2{\rm Re}\{e^{i\delta }\langle Ax_n,x_n\rangle \} <0 $$
as $\xi+\delta\in[\frac{\pi}{2}+\varepsilon,\frac{3\pi}{2}-\varepsilon]$ . { This entails $A\not\perp_w BB^*$. So  $A\not\perp^s_w B$ by Proposition \ref{prop3.1}, a contradiction.}

	(2) If $\dim H=n<\infty$ and $A\in\mathcal{F} $, we divide the proof into two case.
	
{\bf Case 1. } $0\notin {\rm co}\,\{N_A\}$.

	In this case, there exist some $\theta \in[0,2\pi)$ and $\varepsilon >0$ such that
	$$N_A\subseteq \{w(A)e^{i\xi }: \xi \in[\theta +\varepsilon ,\theta +\pi-\varepsilon ]\}.$$ If   $\|B^*x \|\neq 0 $ for all $x\in M_A$, by taking $\delta =\frac{\pi}{2}-\theta$, we have
	 $${\rm Re}\{e^{i\delta }\langle Ax,x\rangle \overline{\langle BB^*x,x\rangle}  \}=\|B^*x \|^2{\rm Re}\{e^{i\delta }\langle Ax,x\rangle \} <0,$$ which implies $A\not\perp_w BB^*$. By Proposition \ref{prop3.1} again,   $A\not\perp^s_w B$, a contradiction. So there exists a vector $x_0\in M_A$ such that $\| B^*x_0\| =0$.
	
{\bf Case 2.}	For any subset $M_A'\subseteq M_A$, there exists a vector  $x_0\in M_A'$ such that $x_0\notin {\rm span}\, \{M_A'\backslash \{x_0\} \}$.

	As $A\perp _w^s B$,   $A\perp _wBC$ holds for all $C\in{\mathcal B}(H)$ by Proposition \ref{prop3.1}. Then,  by \cite[Theorem 3.3]{M},  there exist $t_1,\cdots, t_k>0$, $x_1,\cdots,x_k\in M_A$ ($k\leq 2\dim {\rm span }\,\{BC:C\in{\mathcal B}(H)\}+1$) such that $\sum_{i=1}^{k}t_i=1  $ and
	\begin{equation}\label{3.1}
		\sum_{i=1}^{k}t_i\overline{\langle Ax_i,x_i\rangle} \langle BCx_i,x_i\rangle=0 {\rm \ for \ all \ } C\in{\mathcal B}(H).
	\end{equation}
	Note that, for $\{x_1,\cdots ,x_k\}\subseteq M_A$, there exists a vector $x_i\in \{x_1,\cdots ,x_k\}$ such that $x_i\notin[x_1,\cdots ,x_{i-1},x_{i+1},\cdots ,x_k]$. Now, take $v\in[x_1,\cdots ,x_{i-1},x_{i+1},\cdots ,x_k]^\bot $ and $C=u\otimes v$ with $v\notin[x_i]^\bot $ and $u\in H$.  Eq.(\ref{3.1}) becomes
	$$t_i\overline{\langle Ax_i,x_i\rangle}\langle x_i,v\rangle \langle u,B^*x_i\rangle =0 {\rm \ for \ all \ } u\in H,  $$
	which entails $B^*x_i=0$, as desired.
	
	Combining Case 1 and Case 2, the statement (2) holds.
\end{proof}

Now, by applying Proposition \ref{prop3.3},  we can give another proof of Corollary \ref{prop3.6} for $A,B\in {\mathcal B}_s(H)$.

{\it Proof of Corollary \ref{prop3.6}.} Assume that $A,B\in\mathcal{B}_s (H)$ with $A\not=0$.

(ii)$\Rightarrow $(i): This  is obvious by Proposition \ref{prop3.2}.
	
	(i)$\Rightarrow $(ii):   As $A$ is self-adjoint, $\overline{W(A)} =[a,b]$ for some real numbers $a,b$ with $a\leq b$. We will prove it by considering two cases.
	
{\bf Case 1.}
		$a\neq -b$.

	In this case, we have $A\in\mathcal{E} \cup \mathcal{F} $. So (ii) holds by Proposition \ref{prop3.3}.
	
{\bf Case 2.}
		$a=-b$.

Then $b>0$. Since $A\bot _w^sB$, by Proposition \ref{1.1} and Proposition \ref{prop3.1},  for any $C\in\mathcal{B} (H)$ and any $\theta \in[0,2\pi)$, there exists a sequence $\{x_n^\theta \}\subset H$ such that
	\begin{equation}\label{3.3}
		\lim _{n\rightarrow \infty }|\langle Ax_n^\theta ,x_n^\theta \rangle |=w(A) {\rm \ and \ } \lim _{n\rightarrow \infty }{\rm Re}\{e^{-i\theta }\langle Ax_n^\theta ,x_n^\theta \rangle \overline{\langle BCx_n^\theta ,x_n^\theta \rangle}  \}\geq 0.
	\end{equation}

Still, suppose on the contrary that, for any $\{x_n\}\subset {\rm ball}\,H$ with $\lim\limits _{n\rightarrow \infty }|\langle Ax_n,x_n\rangle |=w(A)=b$, we have $\lim \limits_{n\rightarrow \infty }\|Bx_n \| \neq 0$.
Write the spectral decomposition of $A$ as  $A=\int z{\rm d}E(z)$.   For sufficient small $\varepsilon>0$, let $H_1=E([-b,-b+\varepsilon ])$, $H_3=E([b-\varepsilon,b])$ and $H_2=[H_1\oplus H_3]^\bot $. Then, under the space decomposition $H=H_1\oplus H_2\oplus H_3$,
$$A=A_1\oplus A_2\oplus A_3,\ \ A_i\in\mathcal B (H_i),\ i=1,2,3;$$ and
 $$x_n=\alpha _ny_n\oplus \beta _nz_n\oplus \gamma _nu_n$$ with $y_n\in{\rm ball}\,H_1$, $z_n\in{\rm ball}\,H_2$, $u_n\in{\rm ball}\,H_3$, $\alpha _n,\beta _n,\gamma _n\in\mathbb{C} $ and $|\alpha _n|^2+|\beta _n|^2+|\gamma _n|^2=1$.
By a similar  argument   to that of Case 1 in the proof of Proposition \ref{prop2},  without loss of generality, we may assume that
\begin{equation}\label{3.5}
		\lim _{n\rightarrow \infty }|\gamma  _n|^2=\lim _{n\rightarrow \infty }|\beta  _n|^2=0 {\rm \ and \ }\lim _{n\rightarrow \infty }|\alpha _n|^2=1.
	\end{equation}
Write
	$B=\left(
	\begin{array}{ccc}
		B_1   & B_2   & B_3 \\
		B_2^* & B_4   & B_5 \\
		B_3^* & B_5^* & B_6
	\end{array}
	\right)$.
Take $C=B(-I_1\oplus I_2\oplus I_3)$ and $\theta \in(\frac{\pi}{2}, \frac{3\pi}{2})$. Then
  $$BC=\left(
		\begin{array}{ccc}
			-(B_1^2+B_2B_2^*+B_3B_3^*) & *                       & *                       \\
			*                  & B_2^*B_2+B_4^2+B_5B_5^* & *                       \\
			*                          & *                       & B_3^*B_3+B_5^*B_5+B_6^2
		\end{array}
		\right).$$
    By Eq.(\ref{3.5}), we get
	$$0< \lim _{n\rightarrow \infty }\| Bx_n\|^2=\lim _{n\rightarrow \infty }\langle B^*Bx_n,x_n\rangle
	= \lim _{n\rightarrow \infty }\langle (B_1^2+B_2B_2^*+B_3B_3^*)y_n,y_n\rangle,$$
	which implies
	$$\begin{array}{rl}
		& \lim\limits _{n\rightarrow \infty }{\rm Re}\{e^{-i\theta }\langle Ax_n,x_n\rangle\overline{\langle BCx_n,x_n\rangle}   \}=\lim \limits_{n\rightarrow \infty }{\rm Re}\{e^{-i\theta }b\langle (B_1^2+B_2B_2^*+B_3B_3^*)y_n,y_n\rangle\} \\
		= & b\cos \theta \lim\limits _{n\rightarrow \infty }\langle (B_1^2+B_2B_2^*+B_3B_3^*)y_n,y_n\rangle <0.
	\end{array}$$
	This  contradicts to Eq.(\ref{3.3}). Hence  the statement (ii) holds.
\hfill$\Box$

However, in general,   the converse of Proposition \ref{prop3.2}  is not true  for non normaloid operators $A,B\in\mathcal B(H)$, as shown by (i) of the following result.

\begin{prop}\label{prop3.5}
Assume that $A\in\mathcal{M} _2(\mathbb{C} )\backslash\{0\}$ and  $B\in\mathcal{M} _2(\mathbb{C} )$.

{\rm (i)} If $A^2=0$, that is, $A=\lambda e_1\otimes e_2$  for some $\lambda \in\mathbb{C} $ and $e_1,e_2\in{\rm ball}\,\mathbb{C} ^2$ with $\langle e_1,e_2\rangle=0 $, and if $B=u\otimes v$ is of rank one with $u,v\in\mathbb{C} ^2$ such that  $|\langle u,e_1\rangle |<|\langle u,e_2\rangle |$, then  $A\bot ^s_wB$, but for any $x\in M_A$, $B^*x\not=0$.

{\rm (ii)} If $A^2=0$, that is, $A=\lambda e_1\otimes e_2$  for some $\lambda \in\mathbb{C} $ and $e_1,e_2\in{\rm ball}\,\mathbb{C} ^2$ with $\langle e_1,e_2\rangle=0 $, and if $B=u\otimes v$ is of rank one with $u,v\in\mathbb{C} ^2$ such that  $|\langle u,e_1\rangle |\geq |\langle u,e_2\rangle |$, then $A\bot _w^sB$ if and only if there exists some  $x\in M_A$ such that $B^*x=0$.

{\rm (iii)} If $A^2\not=0$, then $A\bot _w^sB$ if and only if there exists some  $x\in M_A$ such that $B^*x=0$.
\end{prop}

\begin{proof}
By Proposition \ref{prop3.2}, we see that the ``if" parts in the statements (ii) and (iii) are true.

 Also note that both $A\bot _w^sB$ and $A\not=0$ must imply rank$B<2$. Otherwise,   rank$B=2$ and thus is invertible. Then, as   $A\bot _wBC$ for all $C\in\mathcal{M} _2(\mathbb{C} )$ by Proposition \ref{prop3.1}, one gets $A\bot _w^sA$,   contradicting to $A\neq 0$. So we only need to check the statement (i) and the ``only if" parts of the statements (ii) and (iii) for $B$ with rank$B<2$.
	
We will complete the proof by considering two cases.
	
    \begin{case}\label{case3.5.1}
         $A^2=0$. In this case, $A=\lambda e_1\otimes e_2$, where $\lambda \in\mathbb{C} $, $e_1,e_2\in{\rm ball}\,\mathbb{C} ^2$ and $\langle e_1,e_2\rangle=0 $.
    \end{case}

By Proposition \ref{1.1}(1)$\Leftrightarrow $(3) and Proposition \ref{prop3.1}, without loss of generality, we can assume
$A=\left(
  \begin{array}{ccc}
      0 & 1   \\
      0 & 0
    \end{array}
  \right)$
under the orthonormal basis $\{e_1,e_2\}$. It is easy to see that
$$M_A=\{x_\delta :=\frac{1}{\sqrt{2}}
\begin{pmatrix}
  1\\
  e^{i\delta }
  \end{pmatrix}:\delta \in[0,2\pi)\}.$$

{\bf Case 1.1.} $B=u\otimes v$ with $u=\begin{pmatrix}
  u_1\\
  u_2
  \end{pmatrix}$ and  $| u_1|<|u_2|$.

We first  claim that
  \begin{equation}\label{3.2}
  0\in{\rm co}\,\{\overline{\langle Ax_\delta,x_\delta\rangle} \langle (u\otimes q)x_\delta,x_\delta\rangle:\delta \in[0,2\pi)\} \ \ {\rm for \ all}\ \ q=\begin{pmatrix}
  	q_1\\
  	q_2
  \end{pmatrix}\in\mathbb{C}^2. \end{equation}
Note that
$$\overline{\langle Ax_\delta,x_\delta\rangle} \langle (u\otimes q)x_\delta,x_\delta\rangle =\frac{1}{4}(e^{-i\delta }\bar q_1+\bar q_2)(u_1+e^{-i\delta }u_2).$$
If $|q_1|=|q_2|$, then there exists some $\delta \in[0,2\pi)$ such that $e^{-i\delta }\bar q_1+\bar q_2=0$, which means that Eq.(\ref{3.2}) holds.
If $|q_1|\neq |q_2|$, define a function $h$ on \if false $\mathbb{T} =\{e^{i\delta }:\delta \in[0,2\pi)\}$ \fi $\mathbb{D} =\{z\in\mathbb C: |z|\leq 1\}$ by $h(z)=(\bar q_1z+\bar q_2)(u_1+u_2z)$ for all  $z\in\mathbb{D} $. Obviously,  $h$ is analytic on and inside the curve $\mathbb{T} =\{e^{i\delta }:\delta \in[0,2\pi)\}$, and does not vanish on $\mathbb{T}$. Let $N$ stand for the number of zeros of $h$ inside $\mathbb{T}$. Then $N\geq 1$ (as $|u_1|<|u_2|$). According to the argument principle, after the independent variable $z$ makes a complete counterclockwise (positive direction) circuit along $\mathbb{T}$, the change in the argument of the function value $h(z)$ is $2\pi N\geq 2\pi$. So Eq.(\ref{3.2}) also holds.

Now, it follows from  Eq.(\ref{3.2}) that, for any $q\in\mathbb{C}^2 $ and $\theta \in[0,2\pi)$, there exists some $\delta \in[0,2\pi)$ such that ${\rm Re}\,\{e^{i\theta} \overline{\langle Ax_\delta,x_\delta\rangle} \langle (u\otimes q)x_\delta,x_\delta\rangle \}\geq 0$ . By Proposition \ref{1.1}(1)$\Leftrightarrow $(4),   $A\bot_wu\otimes q$ for all $q\in\mathbb{C} ^2$, and so $A\bot _w^su\otimes v$.

However, for any $x\in M_A$, say $x=x_\delta$, we have $B^*x_\delta=\frac{1}{\sqrt{2}}(\bar{u_1}+e^{i\delta}\bar{u_2})v\not=0$.

Therefore, the statement (i) is true.

{\bf Case 1.2.}  $B=u\otimes v$ with $u=\begin{pmatrix}
	u_1\\
	u_2
\end{pmatrix}$ and    $|u_1|\geq |u_2|$.

Assume that $A\bot _w^su\otimes v$. If    $B^*x_\delta=\frac{1}{\sqrt{2}}(\bar u_1+e^{i\delta }\bar u_2)v\neq 0$ for all $\delta \in[0,2\pi)$, then  $|u_1|> |u_2|$.
Also note that, Proposition \ref{prop3.1} implies
$$A\bot _w^su\otimes v\Leftrightarrow A\bot_wu\otimes C^*v {\rm \ for \ all \ }C\in\mathcal{M} _2 (\mathbb{C})\Leftrightarrow A\bot_wu\otimes q {\rm \ for \ all \ }q\in\mathbb{C}^2.$$
 Take
$q_0=\begin{pmatrix}
    0\\
    u_1
    \end{pmatrix}$;
then  $\overline{\langle Ax_\delta,x_\delta\rangle} \langle (u\otimes q_0)x_\delta,x_\delta\rangle=\frac{1}{4}(|u_1|^2+e^{-i\delta }\bar u_1u_2)$. Since $|u_1|>|u_2|$, one gets $|u_1|^2>|u_1u_2|$.  So
$${\rm Re}\,\{e^{i\pi }\overline{\langle Ax_\delta,x_\delta\rangle} \langle (u\otimes q_0)x_\delta,x_\delta\rangle \}={\rm Re}\,\{\frac{1}{4}e^{i\pi }(|u_1|^2+e^{-i\delta }\bar u_1u_2)\}<0 {\rm \ holds\ for \ all\  }\delta \in[0,2\pi),$$
which and  Proposition \ref{1.1}(1)$\Leftrightarrow $(4) gives  $A\not\perp_w u\otimes q_0$,   a contradiction. Hence  there must  exist some $x_\delta \in M_A$ such that $B^*x_\delta=0$.  This completes the proof of the ``only if" part of (ii).

\begin{case}\label{case3.5.1}
  $A^2\not= 0$.
\end{case}

In this case, $W(A)$ is not a closed disk. Assume that $A\bot _w^sB$.  If $0\notin {\rm co}\,\{N _A\}$, then  $A\in\mathcal{F} $. If $0\in {\rm co}\,\{N _A\}$, combining with the fact that $W(A)$ is an ellipse, we conclude that $N _A$ contains exactly two points.  So there exist $x,y\in M_A$ such that $\langle Ax,x\rangle =-\langle Ay,y\rangle$. It follows that $\langle x,y\rangle =0$ by \cite[Problem 1.51]{W2}, which means $A\in\mathcal{F}$. Thus, by Proposition \ref{prop3.3}, there exists some  $x\in M_A$ such that $B^*x=0$,   finishing the proof of ``only if" part of (iii).
\end{proof}

\begin{exmp}  \label{rem2}
Let $A=\left(
  \begin{array}{ccc}
     0 & 1   \\
      0 & 0
    \end{array}
 \right)$ and $B=\left(
 \begin{array}{ccc}
      1 & 1   \\
     2 & 2
   \end{array}
\right)$. Then  $A\bot_w^s B$ by Proposition \ref{prop3.5} (i). However,  $B^*x\neq 0$ for any $x\in M_A$.
\end{exmp}

Recall that $A,B\in\mathcal{B} (H)$ are range orthogonal, denote by $A\bot B$, if $A^*B=0$. It is shown in  \cite{AR3} that,   for any $A,B\in\mathcal{B} (H)$, $A\bot B$ implies $ A\perp ^s_BB $. However, in general, for the strong numerical radius orthogonality, this conclusion no longer holds.

\begin{exmp}
	There exsit matrices $A,B\in\mathcal M_2(\mathbb C)$ such that $A\bot B$, but $A\not\perp ^s_wB$.

In fact, take $A=x\otimes y\in\mathcal{M} _2(\mathbb{C} )$ with $x,y\in{\rm ball}\,\mathbb{C} ^2$ and $|\langle x,y\rangle |\notin\{0,1\}$. Then  $A\in\mathcal{F} $ (\cite[Problem 1.51]{W2}).  Let $B=x^\bot \otimes x^\bot $ with $x^\bot \in{\rm ball}\,\mathbb{C} ^2$ and $\langle x,x^\bot \rangle =0$. It is obvious that $A\bot B$. However, for any unit vector $z\in\mathbb{C} ^2$ with $|\langle Az,z\rangle |=w(A)$, we always have $\| B^*z\| =|\langle z,x^\bot \rangle |\neq 0$, which, by Proposition \ref{prop3.3}, entails $A\not\perp ^s_wB$.
	\end{exmp}

When $A$ is normaloid, the range orthogonality still implies strong numerical radius orthonality.

\begin{prop}\label{prop3.7}

Assume that   $A,B\in\mathcal{B} (H)$ with $A$  normaloid. Then $A\bot B\Rightarrow A\perp ^s_wB.$
\end{prop}

\begin{proof}
    Assume that $A,B\in\mathcal{B} (H)$ with  $A$  normaloid and $A\bot B$. By \cite[Proposition 1]{G}, there exist a sequence $\{x_n\}\subset  {\rm ball}\,H$ and $\lambda \in\mathbb{C} $ such that $|\lambda |=w(A)$ and $\lim\limits _{n\rightarrow \infty } \|(A-\lambda I)x_n \| =0$. Thus, $\lim\limits _{n\rightarrow \infty }|\langle Ax_n,x_n\rangle |=|\lambda |=w(A)$. Note that $B^*A=0$. So
$$0=\lim _{n\rightarrow \infty }\|B^*Ax_n \| =\lim _{n\rightarrow \infty }\|\lambda B^*x_n \|=|\lambda |\lim _{n\rightarrow \infty }\|B^*x_n \|,$$
which means $\lim \limits_{n\rightarrow \infty }\|B^*x_n \|=0$. By Proposition \ref{prop3.2},  $A\bot^s_wB$.
\end{proof}

\section{Characterization of special operators by the strong numerical radius orthogonality}

In this section, we give some characterizations of special bounded linear operators by the strong numerical radius orthogonality.

The first result gives a characterization of the zero operator.

\begin{prop}\label{prop3.9}

 Assume that $A\in\mathcal{B}(H)$. Then  $A=0$ if and only if $X\perp ^s_wA$  for all $X\in\mathcal{B}(H)$.
\end{prop}

\begin{proof}
If $A=0$, then $w(X+AC)\geq w(X)$ holds for all  $X,C\in\mathcal{B} (H)$, and so $X\perp ^s_wA$.

If $X\perp ^s_wA$ holds for all $X\in\mathcal{B}(H)$, then  $A\perp ^s_wA$, that is,  $w(A)\leq w(A+\lambda A)$ for all $\lambda \in\mathbb{C} $. Particulalrly, taking $\lambda =-1$, we obtain $A=0$.
\end{proof}

The following two propositions give  characterizations of the  spectrum of bounded linear operators by the strong numerical radius orthogonality.

\begin{prop}\label{prop3.11}
Assume that $B\in\mathcal{B}(H)$. Then the following statements are equivalent.

   {\rm   (i)}   For any $A\in\mathcal{B} (H)$, we have $A\perp ^s_wB\Rightarrow A=0$.

  {\rm   (ii)}   $0\notin\sigma _{r}(B)$, that is, $B$ is right invertible.
   \end{prop}
   \begin{proof}
   (i)$\Rightarrow $(ii): If $0\in\sigma _{r}(B)$, then  $0\in\sigma _{l}(B^*)$, which implies that  there exists a sequence $\{x_n\}\subset {\rm ball}\,H$ such that $\lim\limits _{n\rightarrow \infty }\|B^*x_n \|=0 $. It follows from Proposition \ref{prop3.2} that $I\bot ^s_wB$,   contradicting to (i). So $0\not\in \sigma _{r}(B)$.

(ii)$\Rightarrow $(i): Assume that  $0\notin\sigma _{r}(B)$. Then $B$ is right invertible. If there exists $A\not=0$ such that $A\bot ^s_wB$, take $C=-C_0A$, where $C_0$ is a right inverse of $B$. It follows that
$0=w(A+BC)\geq w(A)\not=0$, a contradiction. So, we must have $A=0$.
 \if false
We claim that $0\notin \overline{W(BB^*)} $. If not, then there exists a sequence $\{x_n\}\subset {\rm ball}\,H$ such that $\lim\limits _{n\rightarrow \infty }\langle BB^*x_n,x_n\rangle =0$, which means $\lim\limits _{n\rightarrow \infty }\| B^*x_n\| =0$, a contradiction with $0\notin\sigma _{ap}(B^*)$. Since $\sigma (BB^*)\subseteq \overline{W(BB^*)} $, we have $0\notin\sigma (BB^*)$. Take $X\in\mathcal{B} (H)$ such that $BB^*X=I$. For $A\in\mathcal{B} (H)$,  if $A\bot ^s_wB$, then   $A\bot _wBB^*XA$ by Proposition \ref{prop3.1}. So $A\bot _wA$, and hence $A=0$.\fi
   \end{proof}

 \begin{cor}\label{prop3.10}

Assume that $B\in\mathcal{B}_s (H)$.  Then the following statements are equivalent.

  {\rm   (i)}  For any $A\in\mathcal{B} _s(H)$, we have $A\perp ^s_wB\Rightarrow A=0$.

 {\rm   (ii)}    $0\notin\sigma (B)$, that is,  $B$ is invertible.

 \end{cor}

Obviously, the above corollary is also true if both $A$ and $B$ are normal.

\if false
 \begin{proof}
    (i)$\Rightarrow $(ii): Assume that (i) holds, then  $I\not\perp ^s_wB$. So $\lim\limits _{n\rightarrow \infty }\| Bx_n\|\neq 0$ for any $\{x_n\}\subset {\rm ball}\,H$ by Proposition \ref{prop3.6}.
    It follows that $0\notin\sigma (B)$ as $B$ is self-adjoint.

(ii)$\Rightarrow $(i): If $A\perp ^s_wB$ with $0\not=A\in\mathcal{B} _s(H)$, by Corollary \ref{prop3.6} again,  there exists a sequence $\{x_n\}\subset {\rm ball}\,H$ such that $\lim \limits_{n\rightarrow \infty }|\langle Ax_n,x_n\rangle |=w(A)$ and $\lim \limits_{n\rightarrow \infty }\|Bx_n \|=0 $. So $0\in\sigma _{ap}(B)\subseteq \sigma (B)$, a contradiction.
\end{proof}  \fi

Note that the strong numerical radius orthogonality is asymmetric.  So we can introduce the following definition.

   \begin{defn}
   Assume that  $A\in\mathcal{B}(H)$.
We say that  $A$ is   snr-left symmetric if, for $B\in\mathcal{B}(H)$, $A\perp ^s_wB$ implies $B\perp ^s_wA$; and    is   snr-right symmetric if, for $B\in\mathcal{B}(H)$, $B\perp ^s_wA$ implies $A\perp ^s_wB$. Particularly, if $B\in \mathcal{B}_s(H)$, we say that $A$ is $\mathcal{B}_s(H)$-snr-left symmetric and  $\mathcal{B}_s(H)$-snr-right symmetric,   respectively.
 \end{defn}

   The following two results present structures of $\mathcal{B}_s(H)$-snr-left  and  $\mathcal{B}_s(H)$-snr-right symmetric operators.

   \begin{prop}\label{prop3.12}
Assume that $A\in\mathcal{B}_s (H)\backslash \{0\}$.  Then the following statements are equivalent.

  {\rm   (i)}   $A$ is  $\mathcal{B}_s(H)$-snr-left symmetric.

  {\rm   (ii)}   ${\rm rank}\, A=1$.
   \end{prop}
   \begin{proof}
    (ii)$\Rightarrow$(i): Write $A=\lambda x\otimes x$ with $\lambda \in\mathbb{R} \backslash \{0\}$ and $x\in{\rm ball}\,H$. Take any $B\in\mathcal{B} _s(H)$. Assume that $A\bot _w^sB$.  By Corollary \ref{prop3.6}, there exists a sequence $\{x_n\}\subset {\rm ball}\,H$ such that
    \begin{equation}\label{3.6}
    \lim _{n\rightarrow \infty }|\langle Ax_n,x_n\rangle |=w(A) {\rm \ and \ } \lim _{n\rightarrow \infty }\| Bx_n\| =0.
    \end{equation}
Let $x_n=\alpha _nx+\beta _ny_n$, where $y_n\in {\rm ball}([x]^\bot )$, $\alpha _n,\beta _n\in\mathbb{C} $ and $|\alpha _n|^2+|\beta _n|^2=1$. Then $\lim\limits _{n\rightarrow \infty }|\langle Ax_n,x_n\rangle |=\lim\limits _{n\rightarrow \infty }|\alpha _n|^2|\lambda |$, which and Eq.(\ref{3.6}) give $\lim\limits _{n\rightarrow \infty }|\alpha _n|=1$, and so $\lim\limits _{n\rightarrow \infty }|\beta _n|=0$. Note that
$$\| Bx_n\|^2=\langle B^2x_n,x_n\rangle=|\alpha _n|^2\| Bx\| ^2+\alpha _n\bar\beta _n \langle B^2x,y_n\rangle+\bar\alpha _n\beta _n \langle B^2y_n,x\rangle+|\beta _n|^2\langle B^2y_n,y_n\rangle. $$
We infer that $\lim\limits _{n\rightarrow \infty }\| Bx_n\|^2=\| Bx\|^2$. This and Eq.(\ref{3.6}) imply $Bx=0$. So  $B\perp A$.  It follows from  Proposition \ref{prop3.7} that $B\bot _w^sA$.

\if false Let $B=\int \lambda dE(\lambda )$ and $\overline{W(B)}=[\lambda _1,\lambda _2] $. We may assume that $|\lambda _1|\leq |\lambda _2|$ and $\lambda _2>0$. Otherwise, we could replace $B$ by $-B$. Take $\{x_n\}\subseteq {\rm ball}\,H$ with $\lim _{n\rightarrow \infty }|\langle Bx_n,x_n\rangle |=w(B)$.

\begin{case}\label{case3.12.1}
     $|\lambda _1|<|\lambda _2|$.
\end{case}

If $|\lambda _1|<|\lambda _2|$, let $H_1=E[\lambda _2-\varepsilon ,\lambda _2]$ with $\varepsilon >0$, $0\notin[\lambda _2-\varepsilon ,\lambda _2]$ and $|\lambda _1|<|\lambda_2 -\varepsilon |$. With respect to the space decomposition $H=H_1\oplus H_1^\bot $, write $B=B_1\oplus B_2$, $x_n=\alpha _ny_n\oplus \beta _nz_n$ with $y_n\in {\rm ball}\,H_1$, $z_n\in {\rm ball}\,H_1^\bot $, $\alpha _n,\beta _n\in\mathbb{C} $ and $|\alpha _n|^2+|\beta _n|^2=1$. Since
$$|\langle Bx_n,x_n\rangle |=||\alpha _n|^2\langle B_1y_n,y_n\rangle +|\beta _n|^2\langle B_2z_n,z_n\rangle |\leq |\alpha _n|^2\lambda _2+|\beta _n|^2(\lambda _2-\varepsilon ),$$
we arrive at $|\lambda _2|\leq \lambda _2-\varepsilon\lim _{n\rightarrow \infty }|\beta _n|^2 $, which gives $\lim _{n\rightarrow \infty }|\beta _n|^2=0$. Using $0\notin[\lambda _2-\varepsilon ,\lambda _2]$ and $Bx=0$, we infer that $x=0\oplus x_2$ with $x_2\in H_1^\bot $. It follows that
$$\lim _{n\rightarrow \infty }\|Ax_n \| =|\lambda |\lim _{n\rightarrow \infty }|\langle x_n,x\rangle |=|\lambda |\lim _{n\rightarrow \infty }|\beta _n \langle z_n,x_2\rangle|=0,$$
which and Proposition \ref{prop3.6} imply $B\bot ^s_wA$.

\begin{case}\label{case3.12.2}
     $|\lambda _1|=|\lambda _2|$.
\end{case}

Let $H_1=E([-\lambda _2,-\lambda _2+\varepsilon ])$, $H_3=E([\lambda _2-\varepsilon,\lambda _2])$ and $H_2=[H_1\oplus H_3]^\bot $ with $0<\varepsilon <\lambda _2$. For a sequence $\{x_n\}\subseteq {\rm ball}\,H$, write $x_n=\alpha _ny_n\oplus \beta _nz_n\oplus \gamma _nu_n$ with $y_n\in{\rm ball}\,H_1$, $z_n\in{\rm ball}\,H_2$, $u_n\in{\rm ball}\,H_3$, $\alpha _n,\beta _n,\gamma _n\in\mathbb{C} $ and $|\alpha _n|^2+|\beta _n|^2+|\gamma _n|^2=1$. By the same reasoning as in the proof of the Case 2 of Proposition \ref{prop3.6}, we may assume that $\lim _{n\rightarrow \infty }|\alpha   _n|^2=\lim _{n\rightarrow \infty }|\beta  _n|^2=0$. Now to prove $B\bot ^s_wA$, we proceed identically as in the proof of Case 1 of Proposition \ref{prop3.12}???.
\fi

(i)$\Rightarrow$(ii): Assume that for  $B\in\mathcal{B} _s(H)$,  $A\bot _w^sB$ implies $ B\bot _w^sA$.

Suppose on the contrary that ${\rm rank}\,A\geq 2$. Let $A=\int z{\rm d}E(z)$ and $\overline{W(A)}=[a,b] $. Without loss of generality, we may assume that $|a|\leq |b|$ and $b>0$.

\setcounter{case}{0}
\begin{case}\label{case3.12.3}
    $b$ is not an isolated point of $\sigma(A)$.
\end{case}

In this case, take $H_2=E([b-\varepsilon _2,b-\varepsilon _1))$, $H_3=E([b-\varepsilon _1,b])$ and $H_1=(H_2\oplus H_3)^\bot $, where $\sigma (A)\cap [b-\varepsilon _1,b]\neq \emptyset $, $\sigma (A)\cap [b-\varepsilon _2,b-\varepsilon _1)\neq \emptyset $ and $0\notin[b-\varepsilon _2,b]$. Clearly, $H_2,H_3\neq \{0\}$. With respect to the space decomposition $H=H_1\oplus H_2\oplus H_3$, we have $$A=A_1\oplus A_2\oplus A_3\ \ {\rm with}\ \ A_i\in{\mathcal B}(H_i),\ \ i=1,2,3.$$

Take $$B=0\oplus I_2\oplus 0\in{\mathcal B}_s(H).$$ For any $\{x_n\}\subset {\rm ball}\,H$ with $\|x_n\|=1$, write $$x_n=\alpha _ny_n\oplus \beta _nz_n\oplus \gamma _nu_n$$ with $y_n\in{\rm ball}\,H_1$, $z_n\in{\rm ball}\,H_2$, $u_n\in{\rm ball}\,H_3$, $\alpha _n,\beta _n,\gamma _n\in\mathbb{C} $ and $|\alpha _n|^2+|\beta _n|^2+|\gamma _n|^2=1$. On the one hand,  if $\lim \limits_{n\rightarrow \infty }|\langle Ax_n,x_n\rangle |=b$, a direct calculation gives $\lim\limits _{n\rightarrow \infty }\| Bx_n\| =0$. By Corollary \ref{prop3.6}, one has $A\bot ^s_wB$.
On the other hand, if $\lim\limits _{n\rightarrow \infty }|\langle Bx_n,x_n\rangle |=w(B)$, then   $\lim\limits _{n\rightarrow \infty }|\beta _n|=1$ and $\lim\limits _{n\rightarrow \infty }|\alpha  _n|=\lim\limits _{n\rightarrow \infty }|\gamma  _n|=0$. This implies $\lim\limits _{n\rightarrow \infty }\| Ax_n\| =\lim\limits _{n\rightarrow \infty }\| Az_n\|\neq 0$. So $B\not\perp^s_wA$, a contradiciton.

\begin{case}\label{case3.12.4}
     $b$ is an isolated point of $\sigma(A)$ and $\sigma (A)\backslash \{0,b\}\neq \emptyset $.
\end{case}

In this case,   $b\in\sigma _p(A)$ by \cite[Proposition 11.1.4]{C}. Let $H=H_1\oplus H_2\oplus H_3$ with $H_1=E(\{b\})$, $H_3=E(\{0\})$ and $H_2=(H_1\oplus H_3)^\bot $. As $\sigma (A)\backslash \{0,b\}\neq \emptyset $,  $\dim H_2>0$. Under the decomposition $H=H_1\oplus H_2\oplus H_3$,  we have $A=b\oplus A_2\oplus 0$. Now, take $B=0\oplus I_2\oplus 0$. It is easily seen that $A\bot ^s_wB$ and $B\not\perp^s_wA$, a contradiction.

\begin{case}\label{case3.12.5}
    $b$ is an isolated point of $\sigma(A)$ and $\sigma (A)\backslash \{0,b\}= \emptyset $.
\end{case}

Then $\dim \ker({A-bI})\geq 2$. Take any two orthogonal unit vectors $e_1,e_2\in\ker({A-bI})$. Under the decomposition $H=[e_1]\oplus[e_2]\oplus [e_1,e_2]^\bot $, we have $A=b\oplus b\oplus A_3$. By taking $B=1\oplus 0\oplus 0$, one can get  $A\bot ^s_wB$ and $B\not\perp^s_wA$,  a contradiction again.

Now, combining the above three cases, we conclude that rank$A=1$. The proof is finished.
\end{proof}

\begin{prop}\label{prop3.13}
Assume that $A\in\mathcal{B}_s (H)\backslash \{0\}$. Then the following statements are equivalent.

  {\rm   (i)}    $A$ is $\mathcal{B}_s(H)$-snr-right symmetric.

  {\rm   (ii)}    $A$ is invertible.
\end{prop}

\begin{proof}
    (ii)$\Rightarrow $(i): Assume that $A$ is invertible, and  for $B\in\mathcal{B} _s(H)$,  $B\bot ^s_wA$. Then  $0\notin\sigma (A)$, and so  $B=0$   by Corollary \ref{prop3.10}.  It is obvious that  $A\bot ^s_w0$. Hence $A$ is $\mathcal{B}_s(H)$-snr-right symmetric.

    (i)$\Rightarrow $(ii):  Assume that (i) holds. If $A$ is not invertible, then
    $0\in\sigma _{ap}(A)$, which and Corollary \ref{prop3.6} means $I\bot ^s_wA$. On the other hand,  we have  $A\not\perp^s_wI$ by Corollary \ref{prop3.10}. So $A$ is not $\mathcal{B}_s(H)$-snr-right symmetric, a contradiction with (i).
\end{proof}

\if false
\section {Applications: Characterizing maps preserving the (strong) numerical radius orthogonality on operator algebras}

In this section, we give some applications of results obtained in Sections 2-4, which can be appplies to charqacterizing maps preserving the (strong) numerical radius orthogonality on $\mathcal{B} _s(H)$ and $\mathcal{B} (H)$, respectively.

For  maps preserving  numerical radius orthogonality, we have the following results.

\begin{thm}\label{th4.1}
   Assume that $\phi :\mathcal{B} _s(H)\to \mathcal{B} _s(H)$ is an additive surjective map. Then the following statements are equivalent.

(i)  For any $A,B\in\mathcal{B}_s(H)$,  $A\bot _wB$ implies $\phi (A)\bot _w\phi (B)$.

(ii)  For any $A,B\in\mathcal{B}_s(H)$,  $A\bot _{B}B$ implies $\phi (A)\bot _B\phi (B)$.

(iii)  There exist a unitary operator $U\in\mathcal{B} (H)$ and some $\lambda \in\mathbb{R} \backslash \{0\}$ such that $\phi (A)=\lambda UA^\dagger U^*$ for all $A\in\mathcal{B} _s(H)$, where $A^\dagger $ denotes $A$ or $\bar A$.
\end{thm}

\begin{thm}\label{th4.2}
  Assume that $\phi :\mathcal{B}(H)\to \mathcal{B}(H)$ is a linear surjective map. Then for any $A,B\in\mathcal{B}(H)$, $A\bot _wB$ implies $ \phi (A)\bot _w\phi (B)$ if and only if there exist a unitary operator $U\in\mathcal{B} (H)$ and some $\lambda \in\mathbb{C} \backslash \{0\}$ such that $\phi (A)=\lambda UA^\dagger U^*$ for all $A\in\mathcal{B}(H)$, where $A^\dagger $ denotes $A$, $A^t$, $\bar A$ or $A^*$.
\end{thm}


For  maps preserving strong numerical radius orthogonality, we have the following results.

\begin{thm}\label{th4.3}
Assume that $\phi :\mathcal{B} _s(H)\to \mathcal{B} _s(H)$ is a linear surjective map. Then the following statements are equivalent.

(i)  For any $A,B\in\mathcal{B}_s(H)$, $A\bot _w^sB$ if and only if $\phi (A)\bot _w^s\phi (B)$.

(ii)  There exist a unitary opertor $U\in\mathcal{B} (H)$ and $\eta \in\mathbb{R}\backslash \{0\} $ such that $\phi (A)=\eta UA^\dagger U^*$ for all $A\in\mathcal{B} _s(H)$, where $A^\dagger $ denotes  $A$ or $\bar A$.

  Furthermore, if $\dim H<\infty $, then the homogeneity assumption on $\phi $ can be removed.
\end{thm}

\begin{thm}\label{th4.4}
Assume that $\phi :\mathcal{B} (H)\to \mathcal{B} (H)$ is an additive surjective map with $\phi (I)=I$. Then the following statements are equivalent.
	
	(i) For any $A,B\in\mathcal{B} (H)$, $A\bot _w^sB$ if and only if $\phi (A)\bot _w^s\phi (B)$.
	
	(ii)  There exists a unitary operator $U\in\mathcal{B} (H)$ such that  $\phi (A)=UA^\dagger U^*$ for all $A\in\mathcal{B} (H)$, where $A^\dagger $ denotes $A$ or $\bar A$.
\end{thm}

Now, we give the proofs of the above theorems.

{\it Proof of Theorem \ref{th4.1}.}
(i) $\Leftrightarrow $ (ii) is obvious by Proposition \ref{prop3}. (iii) $\Rightarrow  $ (i) is obvious by Proposition \ref{1.1} and the fact $w(A) =w(UAU^*)  =w(U\bar AU^*) $, where $U\in\mathcal{B} (H)$ is a unitary operator.

(ii) $\Rightarrow  $ (iii): Assume that for any $A,B\in\mathcal{B}_s(H)$,  $A\bot _{B}B$ implies $\phi (A)\bot _B\phi (B)$. By \cite[Theorem 3.1]{W}, $\phi (A)=\lambda' \psi (A)$ for all $A\in\mathcal{B} _s(H)$, where $\lambda' \in\mathbb{R} $ and $\psi :{\mathcal B}_s(H)\rightarrow{\mathcal B}_s(H)$ is an isometry. For $\psi$, by \cite[Theorem 2]{K3}, there exist a scalar $\eta \in\{-1,1\}$ and a unitary operator $U\in\mathcal{B} (H)$ such that $\psi (A)=\eta UA^\dagger U^*$ for all $A\in\mathcal{B} _s(H)$, where $A^\dagger $ denotes the $A$ or $\bar A$. Now, take $\lambda =\lambda '\eta $, and then $\lambda \neq 0$ by the surjectivity of $\phi $, as desired.
\hfill$\Box$

{\it Proof of Theorem \ref{th4.2}.}
For the ``if" part, it is clear by the properties of $w(\cdot)$.

For the ``only if" part, since $w(\cdot )$ is a norm of $\mathcal{B} (H)$, by \cite[Theorem 3.1]{BT}, there exist a scalar $\lambda'\in\mathbb{C} $ and a linear isometry $\psi $ on $\mathcal B(H)$ such that $\phi (A)=\lambda'\psi (A)$ for all $A\in\mathcal{B} (H)$. By \cite[Theorem 1]{BH} and the linearity of $\phi $, there exist a complex unit $\mu \in\mathbb{C} $ and a unitary operator $U\in\mathcal{B} (H)$ such that $\psi (A)=\mu UA^\dagger U^*$ for all $A\in\mathcal{B} (H)$, where $A^\dagger $ denotes the $A$, $A^t$, $(A^*)^t$ or $A^*$. Take $\lambda =\lambda '\mu $, we can get $\lambda \neq 0$ by the surjectivity of $\phi $, as desired.
\hfill$\Box$

{\it Proof of Theorem \ref{th4.3}.}
    (ii)$\Rightarrow $(i): Assume that $\phi (A)=\eta UA^\dagger U^*$ for all $A\in\mathcal{B} _s(H)$, where $\eta \in\mathbb{R} \backslash \{0\}$, $U\in\mathcal{B} (H)$ is a unitary operator and $A^\dagger $ denotes the $A$ or $\bar A$. Then for any $A,B\in\mathcal{B} (H)$, we have
$$\begin{array}{rl}
                    & A\perp ^s_wB\Leftrightarrow w(A)\leq w(A+BC) {\rm \ holds \ for \ all \ }C\in\mathcal{B} (H)   \\
    \Leftrightarrow & w(\eta UA^\dagger U^*)\leq w(\eta UA^\dagger U^*+\eta UB^\dagger U^*UC^\dagger U^*){\rm \ holds \ for \ all \ }C\in\mathcal{B} (H) \\
    \Leftrightarrow & w(\phi (A))\leq w(\phi (A)+\phi (B)UC^\dagger U^*){\rm \ holds \ for \ all \ }C\in\mathcal{B} (H)       \\
    \Leftrightarrow & \phi (A)\perp ^s_w\phi (B),
  \end{array}$$
as desired.

(i)$\Rightarrow $(ii): We will prove it by a series of claims.

\begin{note}\label{4.3.1}
     $\phi $ is injective.
\end{note}

Assume that $\phi (X)=0$ for some $X\in\mathcal{B} _s(H)$. As $0\perp ^s_w0$, we have $X\perp ^s_w X$, which means $X=0$. So  $\phi $ is injective.

\begin{note}\label{4.3.2}
     There exist $\xi \in\{-1,1\}$ and an invertible operator $T\in\mathcal{B} (H)$ such that either
     \begin{equation}\label{4.1}
     \phi (X)=\xi TXT^*{\rm \ for \ all \ }X\in\mathcal{B} _s(H);\end{equation}
or
\begin{equation}\label{4.2}
    \phi (X)=\xi TX^tT^*{\rm \ for \ all \ }X\in\mathcal{B} _s(H).\end{equation}
\end{note}

By Proposition \ref{prop3.10} and Claim 1, we infer that
$0\in\sigma (X)\Leftrightarrow 0\in\sigma (\phi (X))$, and so $\phi $ preserves invertibility of  operators.

If $\dim H=\infty $, by \cite[Theorem 6]{LiR}, the claim holds.

If $\dim H<\infty $, as $A-B$ is invertible if and only if $\phi (A-B)=\phi (A)-\phi (B)$ is invertible for all $A,B\in\mathcal{B} _s(H)$, by \cite[Theorem 2.1]{K}, there exist an invertible operator $T\in\mathcal{B} (H)$, $\xi \in\{-1,1\}$ and $C\in\mathcal{B} _s(H)$ such that either
$$\phi (X)=\xi TXT^*+C{\rm \ for \ all \ }X\in\mathcal{B} _s(H);$$
or
$$\phi (X)=\xi TX^tT^*+C{\rm \ for \ all \ }X\in\mathcal{B} _s(H).$$
Note that $\phi (0)=0$. So $C=0$, as desired.

\begin{note}\label{4.3.3}
    The statement (ii) holds, that is, there exist a unitary opertor $U\in\mathcal{B} (H)$ and $\eta \in\mathbb{R}\backslash \{0\} $ such that $\phi (A)=\eta UA^\dagger U^*$ for all $A\in\mathcal{B} _s(H)$, where $A^\dagger $ denotes  $A$ or $\bar A$.
\end{note}

If Eq.(\ref{4.1}) holds, then we have
$$\phi (x\otimes x)=\xi Tx\otimes Tx{\rm \ for \ all \ }x\in{\rm ball}\,H.$$
Since $x\otimes x\bot ^s_wy\otimes y$ for any $y\in[x]^\bot $, then we arrive at $\xi Tx\otimes Tx\bot ^s_w\xi Ty\otimes Ty$. By checking the proof of Proposition \ref{prop3.12}, we have
$$x\otimes x\bot ^s_wB\Leftrightarrow Bx=0$$
for any $x\in H$ and $B\in\mathcal{B} _s(H)$. It follows that $\xi Ty\otimes Ty(Tx)=0$, which means $\langle T^*Tx,y\rangle=0 $. So we have $[x]^\bot \subseteq [T^*Tx]^\bot $ for all $x\in H$, and thus $T^*Tx$ and $x$ are linearly dependent. Then \cite[Lemma 2.2]{H1} applies and tells us that $T^*T=\mu I$ for some $\mu>0$. Let $U=\frac{1}{\sqrt{\mu }}T$ and $\eta =\xi \mu $. Thus, $U$ is a unitary operator and
$$\phi (X)=\xi TXT^*=\xi \mu \frac{1}{\sqrt{\mu }}TX\frac{1}{\sqrt{\mu }}T^*=\eta UXU^* {\rm \ for \ all \ }X\in\mathcal{B} _s(H).$$

Now we assume that Eq.(\ref{4.2}) holds. It follows that
$$\phi (X)=\xi TX^tT^*=\xi TJXJT^*{\rm \ for \ all \ }X\in\mathcal{B} _s(H).$$
Let $\phi '(X)=J\phi (X)J=\xi JTJXJT^*J=\xi T'XT'^*$, where $T'=JTJ\in\mathcal{B} (H)$ is invertible. Note that $A\bot ^s_wB\Leftrightarrow \phi '(A)\bot ^s_w\phi '(B)$ for all $A,B\in\mathcal{B} _s(H)$, a similar reasoning gives that $\phi '(X)=\eta U'XU'^*$ for  all  $X\in\mathcal{B} _s(H)$,
where $\eta \in\mathbb{R} \backslash\{0\}$ and $U'\in\mathcal{B} (H)$ is a unitary operator. Thus,
$$\phi (X)=J\phi '(X)J=\eta JU'XU'^*J=\eta U\bar XU^*$$
for all $X\in\mathcal{B} _s(H)$, where $U=JU'J$ is a unitary operator.

So (ii) holds. The proof of the theorem is completed.
\hfill$\Box$

\begin{proof}
    (ii)$\Rightarrow $(i): Firstly, suppose that $\phi (A)=UAU^*$ for all $A\in\mathcal{B} (H)$, where  $U\in\mathcal{B} (H)$ is a unitary operator. For any $A,B\in\mathcal{B} (H)$, if $A\bot ^s_wB$, we infer that
  $$\begin{array}{rl}
      w(\phi (A))= & w( UAU^*)=w(A)\leq w(  A+ BU^*CU) \\
      = & w( UAU^*+ UBU^*C)=w(\phi (A)+\phi (B)C)
    \end{array}$$
  for all $C\in\mathcal{B} (H)$, which means $\phi (A)\bot ^s_w\phi (B)$.

  Secondly, suppose that $\phi (A)= U\bar AU^*$ for all $A\in\mathcal{B} (H)$, where $U\in \mathcal{B} (H)$ is a unitary operator. For any $A,B\in\mathcal{B} (H)$, if $A\bot ^s_wB$, then we have
  $$\begin{array}{rl}
      w(\phi (A))= & w( U\bar AU^*)=w(A)\leq w( A+ B\overline{U^*CU} )            \\
      =            & w( \bar A+ \bar BU^*CU)=w( U\bar AU^*+U\bar BU^*C)=w(\phi (A)+\phi (B)C)
    \end{array}$$
  for all $C\in\mathcal{B} (H)$, which means $\phi (A)\bot ^s_w\phi (B)$.

For $\phi (A)\bot ^s_w\phi (B)\Rightarrow  A\bot ^s_wB$, the proof is similar and we omit it.

(i)$\Rightarrow $(ii): We will prove it by series of claims.

\begin{note}\label{4.4.1}
     If $\dim H<\infty $, there exist a invertible matrice $A\in\mathcal{B} (H)$ and a automorphism $f$ on $\mathbb{C} $ such that either
     \begin{equation}\label{4.3}
     \phi (T)=Af(T)A^{-1} {\rm \ for \ all \ }T\in\mathcal{B} (H);\end{equation}
or
\begin{equation}\label{4.4}
    \phi (T)=Af(T)^tA^{-1} {\rm \ for \ all \ }T\in\mathcal{B} (H),\end{equation}
where $f(T)$ denotes $[f(t_{ij})]$ for $T=[t_{ij}]$.
\end{note}

Firstly, using Proposition \ref{prop3.9}, we have $X\bot ^s_w0$ for all $X\in\mathcal{B} (H)$. It follows that $\phi (X)\bot ^s_w\phi (0)$ for all $X\in\mathcal{B} (H)$. Applying Proposition \ref{prop3.9} again to $\phi (0)$, we obtain $\phi (0)=0$, which means $\phi $ is a bijection.

Assume that $0\notin\sigma (T^*)$. Note that $\dim H<\infty $, one gets $\sigma (T)=\sigma _p(T)=\sigma _{ap}(T)$ for all $T\in\mathcal{B} (H)$. By Proposition \ref{prop3.11}, we infer that
$$X\bot ^s_wT\Rightarrow X=0 {\rm \ for \ all \ }X\in\mathcal{B} (H),$$
which means $\phi (X)\bot ^s_w\phi (T)\Rightarrow \phi (X)=0$ for all $X\in\mathcal{B} (H)$. Thus, $0\notin\sigma (\phi (T)^*)$ by Proposition \ref{prop3.11} and the surjectivity of $\phi$. So $T$ is invertible if and only if $\phi (T)$ is invertible by the bijectivity of $\phi $. The claim follows immediately by \cite[Corollary 1.3]{FS} and $\phi (I)=I$.

\begin{note}\label{4.4.2}
     If $\dim H<\infty $, the statement (ii) holds.
    \end{note}

We divide the proof into two cases.

\begin{case}\label{case4.4.2.1}
     Eq.(\ref{4.3}) holds.
\end{case}

For any fixed unit vector $x\in H$, there exists some space decomposition $H=[x]\oplus [x]^\bot $ such that $x\otimes x=1\oplus 0$. We obtain $\phi (x\otimes x)=Ax\otimes A^{-1*}x$ by Eq.(\ref{4.3}). Take  $B=0\oplus B_1$, where $B_1$ is invertible, and thus $x\otimes x\bot ^s_wB$, which means $Ax\otimes A^{-1*}x \bot ^s_w Af(B)A^{-1}$. By Proposition \ref{prop3.3} and $Ax\otimes A^{-1*}x\in\mathcal{F} $, there exists a unit vector $x_0\in H$ such that
\begin{equation}\label{4.5}
|\langle (Ax\otimes A^{-1*}x)x_0,x_0\rangle |=w(Ax\otimes A^{-1*}x)\end{equation}
and
\begin{equation}\label{4.6}
(Af(B)A^{-1})^*x_0 =0.
\end{equation}
By Eq.(\ref{4.6}), we have $f(B)^*A^*x_0 =0$. This together with $f(B)^*=0\oplus f(B_1)^*$ entail that $A^*x_0\in[x]$, and thus $x_0=\alpha A^{-1*}x$ for some $\alpha \in\mathbb{C} $. And we have $|\alpha |=\frac{1}{\| A^{-1*}x\| }$. So the Eq.(\ref{4.5}) becomes
$$\begin{array}{rl}
  |\langle \frac{1}{\alpha }Ax,\alpha A^{-1*}x\rangle| &=|\langle (\frac{1}{\alpha }Ax)\otimes (\alpha A^{-1*}x)(\alpha A^{-1*}x),(\alpha A^{-1*}x)\rangle |\\
  &=w((\frac{1}{\alpha }Ax)\otimes (\alpha A^{-1*}x))=\frac{1}{2}(|\langle \frac{1}{\alpha }Ax,\alpha A^{-1*}x\rangle |+\| \frac{1}{\alpha }Ax\| \|\alpha A^{-1*}x \| ),
  \end{array}$$
which means $|\langle \frac{1}{\alpha }Ax,\alpha A^{-1*}x\rangle|=\| \frac{1}{\alpha }Ax\| \|\alpha A^{-1*}x \|$. It follows that $Ax$ and $A^{-1*}x$ are linear dependent for any $x\in H$, which means $A^{-1*}=\lambda A$ for some $\lambda \in\mathbb{R} $. Thus,
$$\phi (T)=Af(T)A^{-1}=\lambda Af(T)A^*=Uf(T)U^* {\rm \ for \ all \  }T\in\mathcal{B} (H),$$
where $U=\sqrt{\lambda }A$ is a unitary operator.

It remains to show that $f$ is a identity or the complex conjugate on $\mathbb{C} $. We may assume that $\phi (T)=f(T)$ for any $T\in\mathcal{B} (H)$. Otherwise, we could replace $\phi (\cdot )$ by $U^*(\cdot)U$. For a fixed orthogonal basis $[e_1,\cdots ,e_n]$ and $\lambda ,\mu \in\mathbb{C} $ with $|\lambda |<|\mu |$, since $\left(
  \begin{array}{ccc}
      \mu & 0       \\
      0   & \lambda
    \end{array}
  \right)\oplus 0\bot ^s_w\left(
  \begin{array}{ccc}
      0 & 0 \\
      0 & 1
    \end{array}
  \right)\oplus I_1$, we infer that
$$\left(
  \begin{array}{ccc}
      f(\mu) & 0           \\
      0      & f(\lambda )
    \end{array}
  \right)\oplus 0\bot ^s_w\left(
  \begin{array}{ccc}
      0 & 0 \\
      0 & 1
    \end{array}
  \right)\oplus I_1.$$
By Proposition \ref{prop3.3}, there exists a unit vector $x_0\in H$ such that
$$|\langle (\left(
    \begin{array}{ccc}
        f(\mu) & 0           \\
        0      & f(\lambda )
      \end{array}
    \right)\oplus 0)x_0,x_0\rangle |=\max\{|f(\mu )|,|f(\lambda )|\}$$
and
$$(\left(
  \begin{array}{ccc}
      0 & 0 \\
      0 & 1
    \end{array}
  \right)\oplus I_1)x_0=0,$$
which means $x_0\in[e_1]$. It follows that $|f(\lambda )|\leq |f(\mu )|$.

Now we assume that $\mu _n\rightarrow 0$. For any $\varepsilon >0$, there exists $\varepsilon _0\in\mathbb{Q} $ and a positive integer $N$ such that $0<\varepsilon _0\leq \varepsilon $ and $|\mu _n|\leq \varepsilon _0$ for any $n>N$. So we have $|f(\mu _n)|\leq f(\varepsilon _0)=\varepsilon _0\leq \varepsilon $, which implies $f(\mu _n)\to 0$. That is, $f$ is a continuous function, which implies that $f$ is a identity or the complex conjugate on $\mathbb{C} $. Thus, the $\phi $ takes the desired forms.

\begin{case}\label{case4.4.2.2}
     Eq.(\ref{4.4}) holds.
\end{case}

If $\phi $ has the form Eq.(\ref{4.4}), by the same reasoning as in the proof above, we can show that
$$\phi (T)=UT^\sharp U^* {\rm \ for \ all \ }T\in\mathcal{M} _n(\mathbb{C} ),$$
where $U$ is a unitary matrix, $T^\sharp $ denotes $T^*$ or $T^t$. We may assume that $U=I$. Otherwise, we can consider $\phi '(\cdot )=U^*(\cdot)U$. For orthogonal unit vectors $x,y\in H$, we have $x\otimes x\bot ^s_w y\otimes x$, but $(x\otimes x)^*\not\perp  ^s_w (y\otimes x)^*$ and  $(x\otimes x)^t\not\perp  ^s_w (y\otimes x)^t$. So the Eq.(\ref{4.4}) can not occur.

\begin{note}\label{4.4.3}
     If $\dim H=\infty $, there exists an invertible bounded linear or conjugate linear operator $A:H\rightarrow H$ such that $\phi (T)=ATA^{-1} {\rm \ for \ all \ }T\in\mathcal{B} (H)$.
\end{note}

By the same reasoning as in the proof of the Claim 1, we have
\begin{equation}\label{4.7}
0\in\sigma _{ap}(T^*)\Leftrightarrow 0\in\sigma _{ap}(\phi (T)^*). \end{equation}
Let $\psi (T)=\phi (T^*)^*$. For any $S\in\mathcal{B} (H)$,  we can write $S=T^*$. We have $0\in\sigma _{ap}(T^*)\Leftrightarrow 0\in\sigma _{ap}(S)$, $0\in\sigma _{ap}(\phi (T)^*)\Leftrightarrow 0\in\sigma _{ap}(\psi (S))$, which and Eq.(\ref{4.7}) imply
$$0\in\sigma _{ap}(S)\Leftrightarrow 0\in\sigma _{ap}(\psi (S)).$$
That is, $S\in\mathcal{B} (H)$ is a left topological divisors of zero if and only if $\psi (S)\in\mathcal{B} (H)$ is a left topological divisors of zero. By \cite[Theorem 3.3 (5)$\Leftrightarrow $(13)]{HC} , there exists an invertible bounded linear or conjugate linear operator $A'$ such that $\psi (T)=A'TA'^{-1}$ for all $T\in\mathcal{B} (H)$. Take $A=A'^{-1*}$, we obtain
$$\phi (T)=\psi (T^*)^*=(A'T^*A'^{-1})^*=A'^{-1*}TA'^*=ATA^{-1} {\rm \ for  \ all \ }T\in\mathcal{B} (H).$$

\begin{note}\label{4.4.4}
     If $\dim H=\infty $, the statement (ii) holds.
\end{note}

We divide the proof into two cases.
\setcounter{case}{0}
\begin{case}\label{case4.4.2.3}
     $A$ is a linear operator.
\end{case}

For any $x\in{\rm ball}\,H$, let $B\in\mathcal{B} (H)$ is satisfied that $\ker B^*=[x]$. Then we have $x\otimes x\bot ^s_wB$, which means $Ax\otimes A^{-1*}x\bot ^s_wABA^{-1}$. Note that $Ax\otimes A^{-1*}x\in\mathcal{E} $, then there exists a sequence $\{x_n\}\subseteq {\rm ball}\,H$ such that
\begin{equation}\label{4.8}
\lim _{n\rightarrow \infty }|\langle (Ax\otimes A^{-1*}x)x_n,x_n\rangle |=w(A)
\end{equation}
and
\begin{equation}\label{4.9}
\lim _{n\rightarrow \infty }\| (ABA^{-1})^*x_n\| =0.
\end{equation}
By Eq.(\ref{4.9}), we have $\lim _{n\rightarrow \infty }\| B^*A^*x_n\| =0$. Let $A^*x_n=\alpha _nx+\beta _ny_n$ with $y_n\in{\rm ball}\, [x]^\bot $, $\alpha _n,\beta _n\in\mathbb{C} $ and $|\alpha _n|^2+|\beta _n|^2=1$. It follows that $B^*A^*x_n=\alpha _nB^*x+\beta _nB^*y_n=\beta _nB^*y_n$, which and $\lim _{n\rightarrow \infty }\| B^*A^*x_n\| =0$ imply $\lim _{n\rightarrow \infty }\beta _n=0$. Note that $\alpha _n=\langle x_n,Ax\rangle $ and ${\rm ball}\, H$ is weak-star compact. So there is a subcolumn, which we still call it $\{x_n\}$, such that $\lim _{n\rightarrow \infty}\langle x_n,Ax\rangle=\lim _{n\rightarrow \infty}\alpha _n=\alpha$ for some $\alpha \in\mathbb{C} $. It follows that $x_n=\alpha _nA^{-1*}x+\beta _nA^{-1*}y_n\rightarrow \alpha A^{-1*}x$, which and $\| x_n\|=1 $ imply $|\alpha |=\frac{1}{\| A^{-1*}x\| }$. Then Eq.(\ref{4.8}) applies and tell us that
$$\lim _{n\rightarrow \infty }|\langle (Ax\otimes A^{-1*}x)x_n,x_n\rangle |=|\langle (Ax\otimes A^{-1*}x)\alpha A^{-1*}x,\alpha A^{-1*}x\rangle |=w(A).$$
By the same reasoning as in the proof of Claim 2, we can show that $Ax$ and $A^{-1*}x$ are linear dependent for any $x\in H$. Thus, $A^{-1*}=\lambda A$ for some $\lambda \in\mathbb{R} $, and then
$$\phi (T)=ATA^{-1}=\lambda ATA^*=UTU^* {\rm \ for \ all \  }T\in\mathcal{B} (H),$$
where $U=\sqrt{\lambda }A$ is a unitary operator, as desired.

\begin{case}\label{case4.4.2.4}
     $A$ is a conjugate linear operator.
\end{case}

Let $\phi '(T)=J\phi (T)J=JATA^{-1}J=BTB^{-1}$, where $B=JA$ is a linear operator. It is obvious that $\phi '$ preserve the strong numerical radius orthogonality in both directions. A similar reasoning gives that $\phi '(T)=U'TU'^* {\rm \ for \ all \ }T\in\mathcal{B} (H)$, where $U'$ is a unitary operator. Then
$$\phi (T)=JU'TU'^{*}J=U\bar TU^*{\rm \ for \ all \ }T\in\mathcal{B} (H),$$
where $U=JU'J$ is a unitary operator.

The proof of the theorem is completed.
\end{proof}

\fi

{\bf Conflict of interest}  The authors declare  that there is no Conflict of interest or Conflict of interest.

\end{document}